\documentclass[reqno]{amsart}
\usepackage{amsmath, amsfonts, epsfig, amssymb, graphics}
\usepackage{hyperref}

{\theoremstyle{plain}
\newtheorem{theorem}{Theorem}[section]
\newtheorem{conjecture}[theorem]{Conjecture}
\newtheorem{proposition}[theorem]{Proposition}
\newtheorem{lemma}[theorem]{Lemma}
\newtheorem{statement}[theorem]{Statement}
\newtheorem{question}[theorem]{Question}
\newtheorem{corollary}[theorem]{Corollary}
}
{\theoremstyle{definition}
\newtheorem{definition}[theorem]{Definition}
\newtheorem{remark}[theorem]{Remark}
\newtheorem{example}[theorem]{Example}
\newtheorem{note}[theorem]{Note}
}

\numberwithin{equation}{section}

\newcommand{\pwp}{$(P,\omega)$-partition}
\newcommand{\si}{\sigma}
\newcommand{\apw}{\mathcal{A}(P,\omega)}
\newcommand{\kpwx}{K_{P,\omega}(x)}
\newcommand{\al}{\alpha}
\newcommand{\be}{\beta}
\newcommand{\lm}{\lambda/\mu}
\newcommand{\plm}{P_{\lm}}
\newcommand{\Or}{O}
\newcommand{\pop}{$(P,\Or)$-partition}
\newcommand{\kpox}{K_{P,\Or}(x)}
\newcommand{\apo}{\mathcal{A}(P,\Or)}
\newcommand{\cuv}{\mathfrak{C}_{vu}}
\newcommand{\mbz}{\mathbb{Z}}
\newcommand{\cc}[2]{\langle #1, #2 \rangle}
\newcommand{\lc}{<_{\mathfrak{C}}}
\newcommand{\C}{C}
\newcommand{\cssf}[1]{s_{#1}}
\newcommand{\la}{\lambda}
\newcommand{\gw}{C}
\newcommand{\gwldmn}{\gw^{\la,d}_{\mu \nu}}
\newcommand{\ldm}{\la/d/\mu}
\newcommand{\inpkn}{\subseteq k \times (n-k)}
\newcommand{\gr}{\mathit{Gr}}
\newcommand{\clmn}{c^\la_{\mu \nu}}
\newcommand{\La}{\Lambda}
\newcommand{\addribbon}[2]{#1[#2]}
\newcommand{\Lam}{\addribbon{\La}{-1}}
\newcommand{\chook}[2]{H_{#2,#1}}
\newcommand{\hnk}{\chook{n-k}{k}}
\newcommand{\D}{D}
\newcommand{\cnk}{\mathfrak{C}_{k,n-k}}
\newcommand{\copr}{\gamma}
\newcommand{\toprib}{R}
\newcommand{\qsym}{\mathit{QSym}}
\newcommand{\vare}{\varepsilon}
\newcommand{\modif}[2]{#1+#2n}
\newcommand{\order}[1]{\langle #1 \rangle}
\newcommand{\Le}{L}
\newcommand{\R}{R}
\newcommand{\IL}{\mathit{IL}}
\newcommand{\SPar}{\mathit{SPar}}
\newcommand{\increm}[3]{#1\uparrow_{#2} #3}
\newcommand{\rb}{\bar{r}}
\newcommand{\dempty}[1]{#1/d/\emptyset}
\newcommand{\error}[2]{B(#1, #2)}

\begin{document}

\title{Cylindric skew Schur functions}

\author{Peter McNamara}
\address{Laboratoire de Combinatoire et d'Informatique Math\'ematique\\
Universit\'e du Qu\'ebec \`a Montr\'eal\\
Case Postale 8888, succursale Centre-ville\\
Montr\'eal (Qu\'ebec) H3C 3P8\\
Canada\\}

\email{mcnamara@lacim.uqam.ca}

\begin{abstract}
Cylindric skew Schur functions, which are a generalisation of skew Schur functions,
arise naturally in the study of $P$-partitions.  Also, recent work of A. Postnikov shows they 
have a strong connection with a problem of considerable current interest: that of
finding a combinatorial
proof of the non-negativity of the 3-point Gromov-Witten invariants.  
After explaining these motivations, we study cylindric skew Schur functions
from the point of view of Schur-positivity.  
Using a result of I. Gessel and C. Krattenthaler,
we generalise a formula of A. Bertram, I. Ciocan-Fontanine and W. Fulton, thus
giving an expansion of an arbitrary cylindric skew Schur function in terms of skew 
Schur functions.  
While we show that no non-trivial cylindric skew Schur functions
are Schur-positive, we conjecture that this can be reconciled using the new concept
of cylindric Schur-positivity.  
\end{abstract}

\maketitle



\section{Introduction}

Cylindric skew Schur functions can be introduced in two very different
ways.  From a combinatorial
perspective, one of these motivations is classical,  
while the other is more contemporary.  The classical motivation begins with 
R. Stanley's
\pwp s, and is centred around
a long-standing conjecture of 
Stanley which gives conditions for a
generating function for the set of \pwp s to be a symmetric function.  
This will be the subject of Section \ref{sec:popartitions} and we will finish the section 
by showing how a natural generalisation of \pwp s leads to the idea
of a cylindric skew Schur function.  We will then give a formal definition of cylindric skew
shapes and cylindric skew Schur functions in Section \ref{sec:cylindricskews}.  At this stage
we make the fundamental observations 
that cylindric skew Schur functions are symmetric functions and that
skew Schur functions are themselves cylindric skew Schur functions.  Therefore, cylindric skew
Schur functions can be viewed as a generalisation of skew Schur functions, and it is
logical to ask which properties of skew Schur functions are preserved under this generalisation.

The contemporary motivation for cylindric skew Schur functions involves the
quantum cohomology ring of the Grassmannian and 
the fundamental open problem of finding a combinatorial proof of the non-negativity of the
3-point Gromov-Witten invariants.  While it will be our starting point in Section
\ref{sec:gwinvariants}, no knowledge of quantum cohomology will be assumed and our
emphasis will be combinatorial.  Gromov-Witten invariants are connected to 
the topic of cylindric skew Schur functions via a theorem of A. Postnikov \cite{Pos04pr}.  Since cylindric
skew Schur functions are symmetric, they can be expanded in terms of Schur functions and 
Postnikov's theorem states that the Gromov-Witten invariants appear as particular coefficients in 
this expansion.  The fundamental open problem mentioned above
then becomes a question about the 
Schur-positivity of cylindric skew Schur functions.  Rather than addressing the open problem
directly, our goal is to give a general study of the Schur-positivity of cylindric skew Schur
functions.  

The geometric definition of Gromov-Witten invariants tells us that cylindric skew Schur
functions in a restricted number of variables are Schur-positive.  In Section 
\ref{sec:schurpositivity}, we show that, except for trivial cases, 
cylindric skew Schur functions are never
Schur-positive in infinitely many variables.  
Since they play an important role in our proof of this result, we investigate the class of 
``cylindric ribbons,'' determining the form of the Schur expansion of their corresponding
cylindric skew Schur functions. 
We also show that, except for trivial cases, cylindric
skew Schur functions are never $F$-positive, where $F$ denotes the fundamental
quasisymmetric functions.  We finish Section \ref{sec:schurpositivity} with a discussion
of the minimum number of variables in which a cylindric skew Schur function will not
be Schur-positive.  

In Section \ref{sec:skewexpansion}, we develop a tool for expanding cylindric 
skew Schur functions as a signed sum of skew Schur functions.  A result of I. Gessel
and C. Krattenthaler \cite{GeKr97} serves as the foundation for our tool, while our
formulation is inspired by a result of A. Bertram, I. Ciocan-Fontanine 
and W. Fulton \cite{BCF99}.

That cylindric skew Schur functions are not Schur-positive is in a sense unfortunate as
we would like an extension of the fact that skew Schur functions are Schur-positive.
In Section \ref{sec:cylschurpos}, we define cylindric Schur-positivity as a natural
generalisation of Schur-positivity and we give evidence in favour of a conjecture that
all cylindric skew Schur functions are cylindric Schur-positive.  

Before beginning in earnest, we introduce terminology and notation that we will
use throughout.  We will denote the sets of integers, non-negative integers
and positive integers by $\mbz$, $\mathbb{N}$ and $\mathbb{P}$
respectively.  We will write $[N]$ to denote the set $\{1,2,\ldots, N\}$.
For symmetric function notation, we will follow \cite{Mac95}.

\begin{figure}
\center
\epsfxsize=25mm
\epsfbox{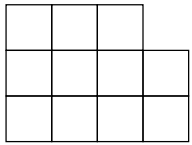}
\caption{The Young diagram for $\la=(4,4,3)$}
\label{fig:youngdiagram} 
\end{figure}

A \emph{composition} of $N$ is a sequence $\al=(\al_1, \ldots, \al_k)$ of 
positive integers that sum to $N$.  
We write $l(\al)$ to denote the 
number of parts of $\al$.
A composition $\la=(\la_1, \ldots, \la_k)$ 
is a \emph{partition} if its sequence of parts is
weakly decreasing.  We also allow partitions to have parts equal to 
zero and we identify $\la$ with the sequence 
$(\la_1, \ldots,\la_k,0,0,\ldots)$.  We write $l(\la)$ for
the number of non-zero parts (length) of $\la$ and $|\la|$ for the sum of 
the parts of $\la$.  
We use $a^k$ in the list of parts of a partition to denote a sequence of $k$
parts of size $a$.  Thus, a partition of the form $(j, 1^k)$ has one part of size $j$
and $k$ parts of size $1$.  Such partitions are called \emph{hooks}.
We let $\emptyset$ denote the unique partition with length 0.
We can represent a partition $\la$ by its Young diagram
drawn in French notation. For example, Figure \ref{fig:youngdiagram}
shows the diagram of the partition $(4,4,3)$. 
The \emph{conjugate partition} $\la'=(\la'_1, \la'_2, \ldots)$ of $\la$
is the partition obtained by reading the column lengths of $\la$, so
in Figure \ref{fig:youngdiagram}, $\la'=(3,3,3,2)$.

If $\mu$ is another partition 
then we say that $\mu \subseteq \la$ if $\mu_i \leq \la_i$
for all $i$.  This is equivalent to saying that the
diagram of $\mu$ is contained in the diagram of $\la$.
If $\mu \subseteq \la$, then
we define the \emph{skew shape}
$\lm$ to be the set of boxes in the
diagram of $\la$ that remain after we remove those boxes corresponding
to the partition $\mu$.  We denote the number of boxes of $\lm$ by $|\lm|$.
A \emph{ribbon} (or \emph{rim hook} or \emph{border strip}) is 
an edgewise connected skew shape that contains
no $2 \times 2$ block of boxes.  An \emph{$n$-ribbon} is then simply a 
ribbon with $n$ boxes.

If a formal power series $f$ can be written uniquely as
a linear combination of
some set of basis elements $\{u_i\}_{i \in I}$ with index set $I$, then we 
write $[u_i]f$ to denote the coefficient of $u_i$ in this linear
combination.  We say that $f$ is \emph{$u$-positive} if
$[u_i]f \geq 0$ for all $i \in I$.  
For example, consider the skew Schur function $s_{\lm}(x)$ in
the variables $x=(x_1, x_2, \ldots)$.  
(For an implicit definition of the skew Schur function $s_{\lm}$, 
see Example \ref{exa:skewshapes}.)  
It can be expanded uniquely in terms of Schur functions
as 
\[
s_{\lm}(x) = \sum_{\nu} \clmn s_\nu (x), 
\]
where $\clmn$ denotes the ubiquitous \emph{Littlewood-Richardson coefficient}.
It is well known that Littlewood-Richardson are non-negative and
skew Schur functions are thus one of the most important examples 
of \emph{Schur-positive} functions.
Schur-positivity has a particular representation-theoretic
significance:
if a homogeneous symmetric function of degree $N$ 
is Schur-positive, then it is known to arise as
the Frobenius image of some representation of the symmetric
group $S_N$.  This is one of the reasons why questions of 
Schur-positivity have received, and continue to receive,
much attention in recent times.  


\section*{Acknowledgements}
This work was begun while the author was a graduate student at MIT.  I am grateful to my
advisor, Richard Stanley, and to Alex Postnikov, for several interesting discussions on 
the topic.  Fran\c{c}ois Bergeron and Christophe
Reutenauer, my mentors at LaCIM, have both made valuable suggestions, and
their expertise and enthusiasm have been of considerable assistance.
Finally, I thank the referee for several suggestions that improved the exposition.


\section{Cylindric skew Schur functions from $(P,\omega)$ partitions}\label{sec:popartitions}

Let $P$ be a finite partially ordered set (poset)
with $N$ elements and let 
\linebreak
$\omega:P \to [N]$
be any bijection labelling the elements of $P$.
We will sometimes refer to elements of $P$ by their images under $\omega$.
The following definition first appeared in \cite{StaThesis}.

\begin{definition}\label{def:pwpartition}
A \pwp\ is a map $\si:P \to \mathbb{P}$ with the following
properties:
\begin{enumerate}
\item[(i)] If $s < t$ in $P$ then $\si(s) \leq \si(t)$. i.e. $\si$ is
order-preserving.
\item[(ii)] If $s < t$ and $\omega(s) > \omega(t)$ then
$\si(s) < \si(t)$.
\end{enumerate}
\end{definition}

Thus a \pwp\ is an order-preserving
map from $P$ to the positive integers
with additional strictness conditions depending on $\omega$.
If $s < t$ is an edge in the Hasse diagram
of $P$ and $\omega(s) > \omega(t)$, then
we will refer to $(s,t)$ as a \emph{strict} edge.  Otherwise,
we will say that $(s,t)$ is a \emph{weak} edge.
In particular, if $\omega$ is itself order-preserving,
then all edges are weak and so
any order-preserving
map from $P$ to $\mathbb{P}$ is a \pwp.
For more information on \pwp s, see \cite{Ges84}, \cite[\S 4.5]{ec1} and
\cite[\S 7.19]{ec2}.
We will denote the set of \pwp s by $\apw$.

Our initial object of study will be the \pwp\ generating function
\linebreak
$\kpwx$ in the variables $x=(x_1, x_2, \ldots)$ defined by
\[
\kpwx = \sum_{\si \in \apw} \prod_{t \in P} x_{\si(t)}
 = \sum_{\si \in \apw} x_1^{\#\si^{-1}(1)}x_2^{\#\si^{-1}(2)}\cdots .
\]
We see that $\kpwx$ is a quasisymmetric function:

\begin{definition}
A \emph{quasisymmetric function} in the variables $x_1, x_2,\ldots$ ,
say with rational coefficients, is a formal power series
$f = f(x) \in \mathbb{Q}[[x_1,x_2,\ldots ]]$ of bounded degree
such that for every
sequence $n_1, n_2, \ldots, n_m \in \mathbb{P}$ of exponents,
\[
\left[x_{i_1}^{n_1} x_{i_2}^{n_2} \cdots x_{i_m}^{n_m}\right]f = 
\left[x_{j_1}^{n_1} x_{j_2}^{n_2} \cdots x_{j_m}^{n_m}\right]f
\]
whenever $i_1 < i_2 < \cdots <i_m$ and
$j_1 < j_2 < \cdots <j_m$.
\end{definition}

Notice that we get the definition of a symmetric function when we change
the condition that the sequences $i_1, i_2, \ldots,i_m$
and $j_1, j_2, \ldots, j_m$ be strictly increasing
to the weaker condition that each sequence consists of distinct elements.
As an example, the formal power series
\[
        f(x) = \sum_{1\leq i < j} x_i^2 x_j
\]
is a quasisymmetric function but is not a symmetric function.
While they appeared implicitly in earlier work, quasisymmetric 
functions were first defined by Gessel \cite{Ges84}, motivated
by the function $\kpwx$.  

There are two bases for quasisymmetric functions that will
be useful.
If 
\linebreak
$\al=(\al_1, \ldots , \al_k)$ is a composition of $N$,
then we define the
\emph{monomial quasisymmetric function} $M_{\al}$ by
\begin{equation}
M_{\al} = \sum_{1\leq i_1 < \cdots <i_k} 
                x_{i_1}^{\al_1} \cdots x_{i_k}^{\al_k} .
\end{equation}
It is clear that the set $\{M_\al\}$, where $\al$ ranges over all
compositions of $N$, forms a basis for the vector space of
quasisymmetric functions that are homogeneous 
of degree $n$. 
We also define
the \emph{fundamental quasisymmetric function} $F_\al$ by
\[
F_\al = \sum_{\be \geq \al} M_\be,
\]
where $\be \geq \al$ denotes that $\be$ is a composition
of $N$ that is a refinement of the composition $\al$.  For example,
$F_{31} = M_{31} + M_{211} + M_{121} + M_{1111}$.
By Inclusion-Exclusion, 
\[
M_\al = \sum_{\be \geq \al} (-1)^{l(\be)-l(\al)} F_\be .
\]
Hence the set $\{F_\al\}$, where $\al$ ranges over all 
compositions of $N$, forms an alternative basis for the 
vector space of homogeneous quasisymmetric functions of degree $N$.

We are now ready to give a concrete example of $\kpwx$.

\begin{example}\label{exa:pwpartition}
Suppose $(P,\omega)$ is given by Figure \ref{fig:pwpartition},
where the double edges correspond to strict edges of $P$.
\begin{figure}
\center
\epsfxsize=25mm
\epsfbox{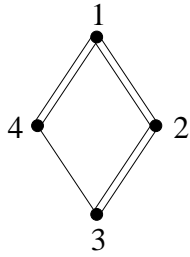}
\caption{A poset $P$ with its labelling $\omega$}
\label{fig:pwpartition} 
\end{figure}
\begin{table} 
\center
\begin{tabular}{cc}
Values of $\si$ & Contribution to $\kpwx$ \\
\hline
$\si(3) = \si(4) < \si(2) < \si(1)$ & $M_{211}$ \\
$\si(3) < \si(4) = \si(2) < \si(1)$ & $M_{121}$ \\
$\si(3) < \si(4) < \si(2) < \si(1)$ & $M_{1111}$ \\
$\si(3) < \si(2) < \si(4) < \si(1)$ & $M_{1111}$
\end{tabular}
\caption{\pwp s for Figure \ref{fig:pwpartition}}
\label{tab:relativeorderings}
\end{table}
We see that a \pwp\ $\si$ must fall into exactly one of the classes shown
in Table \ref{tab:relativeorderings}.  
We conclude that 
\[
\kpwx = M_{211} + M_{121} + 2M_{1111} = F_{211} + F_{121}.
\]
Therefore, the monomial $x_1^2 x_2 x_3$ appears with coefficient 1 in
$\kpwx$ whereas $x_1 x_2 x_3^2$ has coefficient 0.
In particular, $\kpwx$ is not symmetric.
In general, suppose that we have a quasisymmetric function
$f = \sum c_\al M_\al$, where the sum is over all compositions
$\al$ of $N \in \mathbb{P}$.
We see that $f$ is a symmetric function if and only if
$c_\al = c_\be$ whenever $\al$ and $\be$ are compositions
with the same multiset of parts.
\end{example}

\begin{example}\label{exa:skewshapes}
Let $\lm$ be a skew shape with $|\lm|=N$.
We define a \emph{Schur labelling} of $\lm$ to be a labelling 
of the boxes of $\lm$ with the numbers $[N]$ that increases
down columns and from left to right along rows.
Given a Schur labelling $\omega$ of 
$\lm$, let $(\plm, \omega)$ denote the labelled poset
suggested by rotating the boxes of $\lm$ by $45^\circ$ counterclockwise.
These definitions are best explained by an example and Figure
\ref{fig:skewshape} shows a Schur labelling $\omega$ of $\lm$ and
the corresponding labelled poset $(\plm, \omega)$.
\begin{figure}
\center
\epsfxsize=80mm
\epsfbox{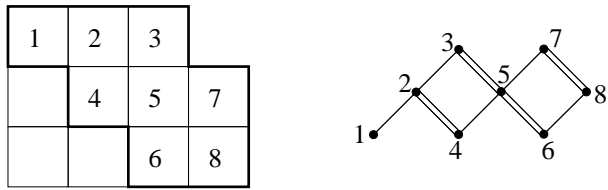} 
\caption{A Schur labelled skew shape and its corresponding labelled poset}
\label{fig:skewshape}
\end{figure}
We say that $(\plm,\omega)$ is a \emph{Schur labelled skew shape poset}
or just a \emph{skew shape poset}.

We see that a \pwp\ of a skew shape poset $(\plm,\omega)$ corresponds
to an assignment of positive integers to the boxes of
$\lm$ that weakly increases from left to right along rows and
strictly increases up columns.
This is exactly the definition of a
\emph{semistandard Young tableau}
of shape $\lm$.
Therefore, the quasisymmetric function $K_{\plm,\omega}$ gives us
exactly the \emph{Schur function} $s_{\lm}$.
We conclude that $\kpwx$ is symmetric
if $(P,\omega)$ is a skew shape
poset.
\end{example}

This brings us to Stanley's $P$-partitions Conjecture \cite{StaThesis}.
We say that two labelled posets are isomorphic if there exists a poset 
isomorphism between them that sends weak edges to weak edges
and strict edges to strict edges.

\begin{conjecture}\label{con:stanleysppartitions}
Let $(P,\omega)$ be a labelled poset.
$\kpwx$ is symmetric if and only if $(P,\omega)$ is isomorphic to a
Schur labelled skew shape poset.
\end{conjecture}

In \cite[Exercise 4.23]{ec1} and \cite{ec1errata}, this conjecture
is shown to be true when $\omega$ is a linear extension.
Using \cite{SteSoftware}, we have verified
the conjecture for all posets $P$ with $|P| \leq 8$.

The reader may already have observed that to calculate $\kpwx$,
we don't need to know the full labelling $\omega$.  It suffices
to know which edges are strict and which edges are weak.
Therefore, from now on, we will often omit the labels on the vertices,
and when we refer to 
a ``labelled poset,''
we mean a poset with strict and weak edges which
can come from some underlying labelling.

\begin{figure}
\center
\epsfxsize=80mm
\epsfbox{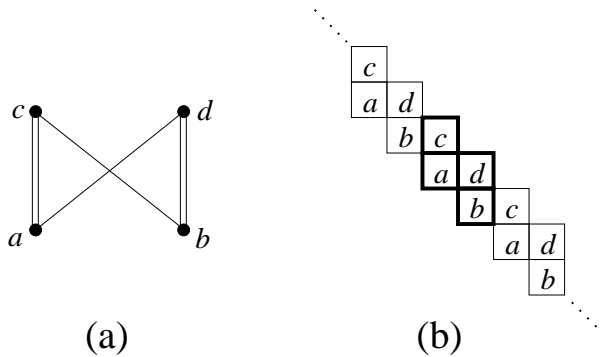}
\caption{An oriented poset and its corresponding cylindric skew shape}
\label{fig:cylindric}
\end{figure}

However, as we shall see, not all designations of strict
and weak edges can come from a labelling. Really though, 
it seems natural
that given a poset $P$, we should allow ourselves to choose
any designation $\Or$ of strict and weak edges.  We will then refer to 
$(P,\Or)$ as an \emph{oriented poset}, with labelled posets
themselves considered to be a special class of oriented posets.  
For example, consider the oriented poset $(P, \Or)$ shown in 
Figure \ref{fig:cylindric}(a) and suppose that it actually corresponds
to a labelled poset $(P, \omega)$.  
Then $\omega$ would have to satisfy $\omega(a) > \omega(c) > \omega(b)
> \omega(d) > \omega(a)$, which is impossible.
With this example in mind, given an oriented poset $(P,\Or)$, 
suppose we think of the Hasse diagram of $P$
as a directed graph, with strict edges oriented upwards, and weak edges
oriented downwards.
We then define a \emph{cycle} of the oriented poset
$(P,\Or)$ to be a cycle in the Hasse diagram of $P$ viewed in this
way as a directed graph.  So a cycle in an oriented poset can be thought of
as a closed path that ``goes up'' on strict edges and down on weak
edges.  Note that in a labelled poset, edges will always be oriented towards
the smaller label.  It follows that if $(P,\Or)$ is a labelled poset, then
it has no cycles.  Furthermore, the converse can also
be shown to be true: an oriented poset is a labelled poset if it has no cycles.  

We define a \pop\ in the obvious manner:

\begin{definition}
Let $P$ be a poset with a designation $\Or$ of strict and weak edges.
A \pop\ is a map $\si:P \to \mathbb{P}$ with the following
properties:
\begin{enumerate}
\item[(i)] If $s < t$ in $P$ then $\si(s) \leq \si(t)$. i.e. $\si$ is
order-preserving.
\item[(ii)] If $s < t$ and $(s,t)$ is a strict edge, then
$\si(s) < \si(t)$.
\end{enumerate}
\end{definition}

As one might expect, we denote the set of
\pop s by $\apo$ and 
we define the generating function $\kpox$ analogously
to $\kpwx$:
\[ 
\kpox = \sum_{\si \in \apo} x_1^{\#\si^{-1}(1)}x_2^{\#\si^{-1}(2)}\cdots .
\]
Isomorphism of oriented posets is defined exactly as for labelled posets.

\begin{example}\label{exa:bowtie}
Consider again the oriented poset $(P,\Or)$
shown in Figure \ref{fig:cylindric}(a).
Using the same method as in Example \ref{exa:pwpartition}, we can compute
that 
\begin{eqnarray*}
\kpox & = & M_{22} + 2M_{211} + 2M_{121} + 2M_{112} + 4M_{1111} \\
& = & F_{22} + F_{211} + 2F_{121} + F_{112} - F_{1111}.
\end{eqnarray*}
Since $(P,\Or)$ is not a labelled poset, it is not considered in 
Conjecture \ref{con:stanleysppartitions}.  Even so, since $\kpox$ 
is symmetric, we might wonder if it somehow comes from a skew shape.
Referring now to Figure \ref{fig:cylindric}(b), we see 
that we need the box corresponding to $a$ to be
directly below the box corresponding to $c$ and directly to the left
of the box corresponding to $d$.  Also, we need the box corresponding
to $b$ to be directly below the box corresponding to $d$ and directly
to the left of the box corresponding to $c$.  Naively putting this
all together, we might be led to the construction
in Figure \ref{fig:cylindric}(b).  We refer to such constructions as
\emph{cylindric skew shapes} and then $\kpox$ is an example of
a \emph{cylindric skew Schur function}.  
This example motivates the
formal definitions of the next section.
\end{example}


\section{Cylindric skew Schur functions}\label{sec:cylindricskews}

Cylindric skew shapes are not a new idea and there are three references in particular
that are of great relevance to our work.  The first of these is 
\cite{GeKr97}, which will play an important role in Section \ref{sec:skewexpansion}.
Semistandard cylindric tableaux, which we will shortly define, appear under the
name ``proper tableaux'' in \cite{BCF99}.  The main result of \cite{Pos04pr}
serves as the starting point for our results of the next section.  
Also, for the following introduction to the 
notation and definition of cylindric skew shapes, we will largely follow \cite{Pos04pr}.

Fix positive integers $u$ and $v$.  We define the
\emph{cylinder}
$\cuv$ to be the following quotient of the integer lattice $\mbz^2$:
\[
\cuv = \mbz^2 / (-u, v)\mbz .
\]
In other words, $\cuv$ is the quotient of 
$\mbz^2$
modulo a shifting action which sends $(i,j)$ to $(i-u,j+v)$.
For $(i,j) \in \mbz^2$,  we let 
$\cc{i}{j} = (i,j) + (-u,v)\mbz$
denote the corresponding element of $\cuv$.
$\cuv$ inherits a natural partial order $\leq_{\mathfrak{C}}$ from $\mbz^2$
which is generated by the relations $\cc{i}{j} \lc \cc{i+1}{j}$ and
$\cc{i}{j} \lc \cc{i}{j+1}$.

Note that this partial order is antisymmetric since $u$ and $v$ are positive.
Recall that a subposet $Q$ of a poset $P$ is said to be \emph{convex} if, for all elements
$x < y < z$ in $P$, we have $y \in Q$ whenever we have $x, z \in Q$.  

\begin{definition}\index{cylindric skew shape}
A \emph{cylindric skew shape} is a finite convex subposet of
the poset $\cuv$.
\end{definition}

\begin{example}\label{exa:skewshapecylindric}
We can regard skew shapes $\lm$ as a special case of cylindric skew shapes.
Suppose $\lm$ fits inside a box of height $v$ and width $u$.
We embed $\lm$ in $\cuv$ by mapping the box in the $i$th row and $j$th
column of $\lm$ to $\cc{i}{j}$. 
Figure \ref{fig:skewshapecylindric} shows the resulting image of
$\lm$ in $\mbz^2$, with one representative of $\lm$ shown in bold.
Notice that elements of different representatives
of $\lm$ are always incomparable in $\mbz^2$.
Of course, we could also embed $\lm$ in $\mathfrak{C}_{v'u'}$
where $v' \geq v$ and $u' \geq u$. 
\begin{figure}
\center
\epsfxsize=75mm
\epsfbox{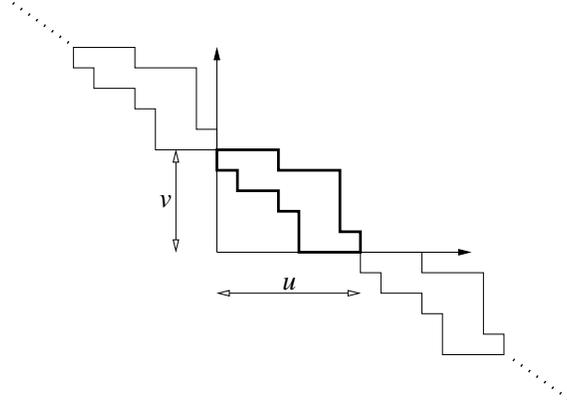}
\caption{Demonstration that skew shapes are cylindric skew shapes}
\label{fig:skewshapecylindric}
\end{figure}
\end{example}

\begin{example}
The class of \emph{cylindric ribbons} will play an important role and they are defined 
in the analogous way to ribbons in the classical case.  
As we just did for skew shapes, we will identify any cylindric
skew shape with its corresponding set of boxes in $\mbz^2$. 
Note that the 
skew shapes from the previous example can be edgewise connected when viewed as subsets
of $\cuv$.  However, they are not edgewise connected when viewed as subsets of $\mbz^2$, 
as in the figure.  

\begin{definition}
A \emph{cylindric ribbon} is a cylindric skew shape which, when viewed as a subset
$\mbz^2$, is edgewise connected and contains no $2 \times 2$ block of boxes.
\end{definition}
\noindent The cylindric skew shape in Figure \ref{fig:cylindric}(b) is an example of a cylindric ribbon.  
\end{example}

Suppose $\C$ is a cylindric skew shape which is a subposet
of the cylinder $\cuv$.
Let us define what we mean by the rows and columns of $\C$.
The \emph{$p$-th row} is the set $\{\cc{i}{j} \in \C \ |\ j=p\}$
and the \emph{$q$-th column}
is the set $\{\cc{i}{j} \in \C \ |\ i=q\}$.
\footnote{
In \cite{Pos04pr}, rows and columns are defined the other way around.}
So the rows only depend on $p \bmod v$ and the columns only
depend on $q \bmod u$.  Thus the cylinder $\cuv$ has exactly
$v$ rows and $u$ columns.

\begin{definition}
For a cylindric skew shape $\C$, a 
\emph{semistandard cylindric tableau} of shape $\C$ is a map $T: \C \to \mathbb{P}$ that
weakly increases in the rows of $\C$ and strictly increases in the columns.
\end{definition}

\begin{figure}
\center
\epsfxsize=110mm
\epsfbox{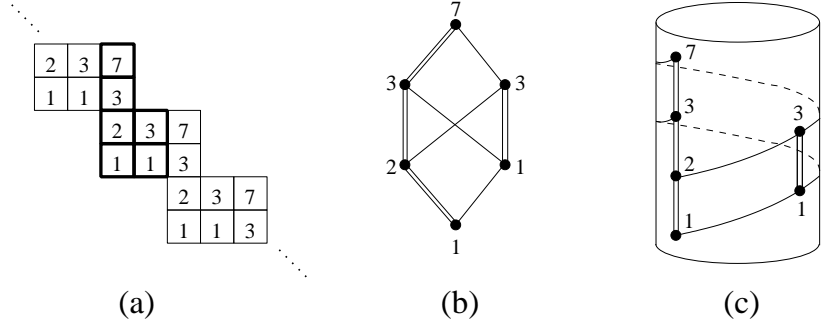}
\caption{A cylindric tableau and its cylindric skew shape poset}
\label{fig:cylindrictableau}
\end{figure}

See Figure \ref{fig:cylindrictableau}(a) for an example.
We are now ready to define our main object of study.

\begin{definition}
For a cylindric skew shape $\C$, the \emph{cylindric skew Schur function} $\cssf{\C}(x)$
in the variables $x=(x_1, x_2, \ldots )$ is defined by
\[
\cssf{\C}(x) = \sum_T \prod_{c \in \C} x_{T(c)}
 = \sum_T x_1^{\#T^{-1}(1)}x_2^{\#T^{-1}(2)}\cdots ,
\]
where the sums are over all semistandard cylindric tableaux $T$ of shape $\C$.
\end{definition}

The terminology ``cylindric skew Schur function'' is partially justified by the following
two observations:

\begin{example}
Because of Example \ref{exa:skewshapecylindric}, skew Schur functions and, in particular, 
Schur functions are all examples of cylindric skew Schur functions.
\end{example}

\begin{theorem}
For any cylindric skew shape $\C$, $\cssf{\C}(x)$ is a symmetric function.
\end{theorem}

We omit the proof as it is basically the same as the proof from \cite{BeKn72}, which also
appears as \cite[Theorem 7.10.2]{ec2}, 
that the skew Schur function $s_{\lm}(x)$ is symmetric.

When $\C$ is a cylindric skew shape which is a subposet of the cylinder $\cuv$ with $u,v \geq 2$, 
we can give a definition of $\cssf{\C}(x)$ in terms of \pop s.  
Now the elements of $\C$ inherit a partial order 
from $\cuv$.
Suppose we consider the vertical edges $\cc{i}{j} \lc \cc{i}{j+1}$
of $\C$
to be strict and the horizontal edges $\cc{i}{j} \lc \cc{i+1}{j}$
to be weak.  This designation of strict and weak edges makes $\C$ into
an oriented poset $(P,\Or)$, which we refer to as a
$\emph{cylindric skew shape poset}$.
  We see that the generating function 
$\kpox$ then coincides exactly with $\cssf{\C}(x)$.
We will find it convenient to switch to this viewpoint of $\C$ and $\cssf{\C}(x)$ at times. 

For example, we encountered a cylindric skew shape poset in
Figure \ref{fig:cylindric}.
Also, because of Example \ref{exa:skewshapecylindric}, skew
shape posets
are always cylindric skew shape posets.
As a further example, 
Figure \ref{fig:cylindrictableau} shows a semistandard cylindric tableau as
well as the corresponding cylindric skew shape poset $(P, \Or)$,
with elements labelled by their images under the corresponding \pop.
Figures 5(b) and 5(c) show the same poset, but the intention of Figure 5(c) is to 
justify the use of the word ``cylindric.''

\begin{note}\label{not:uvequal1}
In the definition above of cylindric skew shape posets,
we required that $u, v \geq 2$.  This is to ensure that
$\cc{i}{j} \lc \cc{i+1}{j}$ and $\cc{i}{j} \lc \cc{i}{j+1}$ are
covering relations.  Indeed, suppose that $u=1$ and $v>1$.  Then we would
have
\[
\cc{0}{0} \lc \cc{0}{1} \lc \cdots \lc \cc{0}{v} = \cc{1}{0}
\]
and so $\cc{0}{0}$ is not covered by $\cc{1}{0}$.  We have a similar
problem if $v=1$.  We will occasionally have a need to consider the
cases when $u$ or $v$ is 1, but we will deal with these cases separately.   
\end{note}

We wish to conclude this section by mentioning some computations with
oriented posets that might affect one's belief in the truth of Conjecture
\ref{con:stanleysppartitions}.
Based on their construction, one could argue that cylindric skew shape posets
play the same role for oriented posets as skew shape posets do for 
labelled posets.  In fact, the following two theorems make this analogy 
even more concrete.

\begin{figure}
\center 
\epsfxsize=90mm
\epsfbox{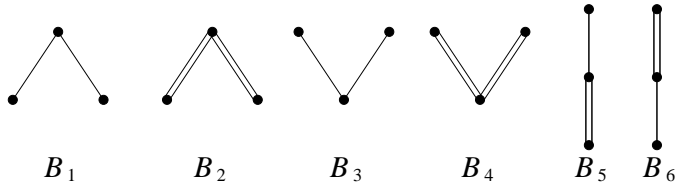} 
\caption{The six ``forbidden'' convex subposets}
\label{fig:triplets} 
\end{figure}

The requirement that $(P,\omega)$ be a skew shape poset seems, in effect, 
to be a  global condition on $(P, \omega)$.  The following result of C. Malvenuto
shows that being a skew shape poset can, in fact, be expressed as a local
condition.  The proof follows from \cite{Mal93}, with some clarification 
and further analysis of her results in \cite{MalThesis}.  Consider the six 3-element
posets $B_1, B_2, \ldots, B_6$ shown in Figure \ref{fig:triplets}.

\begin{theorem}
Let $(P,\omega)$ be a labelled poset.  Then $(P,\omega)$ is isomorphic to
a skew shape poset if and only if $(P,\omega)$ does not contain any $B_i$ as
a convex subposet.
\end{theorem}

It follows that proving Conjecture \ref{con:stanleysppartitions} boils down to 
showing that if $(P,\omega)$ contains a $B_i$, then $\kpwx$ is not
symmetric.
We now state the analogous result for oriented posets.  The proof uses several
of Malvenuto's ideas as well as some new ones, and can be found in 
\cite{McNThesis}.

\begin{theorem}\label{thm:malvenutoextended}
Let $(P,\Or)$ be an oriented poset.  Then every connected component of $(P,\Or)$
is isomorphic to a cylindric skew shape poset if and only if $(P,\Or)$ does not
contain any $B_i$ as a convex subposet.
\end{theorem}

Based on the similarity of these two theorems and other evidence, it is natural
to think that the following analogy of Conjecture \ref{con:stanleysppartitions} might be true:

\begin{statement}\label{sta:falseconjecture}
Let $(P,\Or)$ be an oriented poset.  $\kpox$ is symmetric if and only if every 
connected component of $(P,\Or)$ is isomorphic to a cylindric skew shape
poset.
\end{statement}

\begin{figure}
\center 
\epsfxsize=90mm
\epsfbox{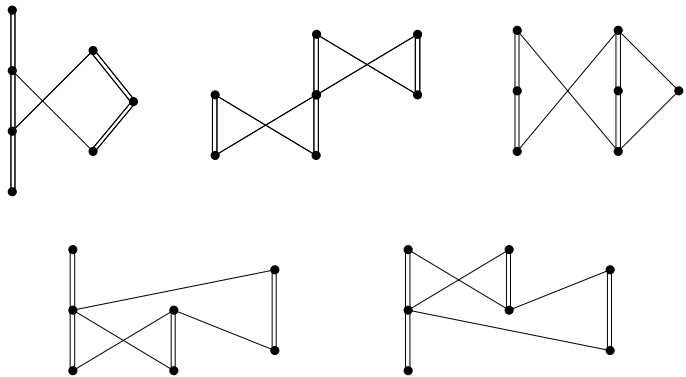} 
\caption{Counterexamples to Statement \ref{sta:falseconjecture}}
\label{fig:counterexamples} 
\end{figure}

Because of Theorem \ref{thm:malvenutoextended}, proving this statement
boils down to showing that if $(P,\Or)$ contains a $B_i$, then 
$\kpox$ is not symmetric, just like for Conjecture \ref{con:stanleysppartitions}.
However, Statement \ref{sta:falseconjecture} is false.  The smallest counterexamples
have 7 elements, and are shown in Figure \ref{fig:counterexamples}.  They
were found using \cite{SteSoftware}.

This might cause one to question the validity of Conjecture \ref{con:stanleysppartitions}. 
On the other hand, there are other things that are true for labelled posets $(P,\omega)$
but not for general oriented posets.  For example, $\kpwx$ has a nice expansion,
with all non-negative integer coefficients, in terms of the basis of 
fundamental quasisymmetric functions $F_\al$.
(See \cite[Corollary 7.19.5]{ec2}.)


\section{Cylindric skew Schur functions from Gromov-Witten invariants}\label{sec:gwinvariants}

As mentioned previously, there is an entirely different -- and relatively new --
reason to be interested in cylindric
skew Schur functions.  This motivation is centred around the main result of
\cite{Pos04pr}.  A nice introduction,
with emphasis on the context and the importance of Postnikov's result can be found in 
\cite{Sta02}.  Here, however, we merely extract from these two references
the minimum amount of background
necessary to show how Postnikov's work ties together cylindric skew Schur functions
and an important open problem in Quantum Schubert Calculus.  

Given $k$ and $n$ with $n > k \geq 1$, we let $\gr_{kn}$ denote the manifold
of $k$-dimensional subspaces of $\mathbb{C}^n$.  $\gr_{kn}$ is a complex
projective variety known as the \emph{Grassmann variety} or \emph{Grassmannian}.  
For a partition $\la$, we will write $\la \inpkn$ if the Young diagram 
for $\la$ has at most $k$ rows and at most $n-k$ columns.
In this case, we let $\la^\vee$ denote the partition $(n-k-\la_k, \ldots, n-k-\la_1)$.
Given $\la, \mu, \nu \inpkn$, we let $\gwldmn$ denote the 
\emph{(3-point) Gromov-Witten invariant}, defined geometrically as the number 
of rational curves of degree $d$ in $\gr_{kn}$ that meet fixed generic translates
of the Schubert varieties $\Omega_{\la^\vee}, \Omega_\mu$ and $\Omega_\nu$,
provided that this number is finite.  This last condition implies that $\gwldmn$ is
defined if $|\mu| + |\nu| = nd + |\la|$, and otherwise we set $\gwldmn=0$.
If $d=0$, then a degree 0 curve is just a point in $\gr_{kn}$ and we get
the geometric interpretation of the Littlewood-Richardson coefficient
$\clmn = \gw^{\la, 0}_{\mu \nu}$.  While we do not claim that this paragraph
is sufficient to give a firm understanding of $\gwldmn$, we do claim that it
is clear from this geometric definition that $\gwldmn \geq 0$.  
No algebraic or combinatorial proof of this inequality is known and, as 
stated in \cite{Sta02}, it 
is a fundamental open problem to find such a proof.

Postnikov's result shows that the Gromov-Witten invariants $\gwldmn$ appear as
the coefficients when we expand certain cylindric skew Schur functions in terms
of Schur functions.  It follows that improving our understanding of this expansion
could lead to a solution of the open problem. 

Before stating his result, we need to introduce some notation
that will allow us to write any cylindric skew shape in the form 
$\ldm$, where $\la$ and $\mu$ are partitions and where $d \in \mathbb{N}$.
From this point on, unless otherwise stated, all of our cylindric skew shapes $\C$ 
will be subposets of the cylinder $\cuv$ with $v= k$ and $u = n-k$.

\begin{figure}
\center 
\epsfxsize=90mm
\epsfbox{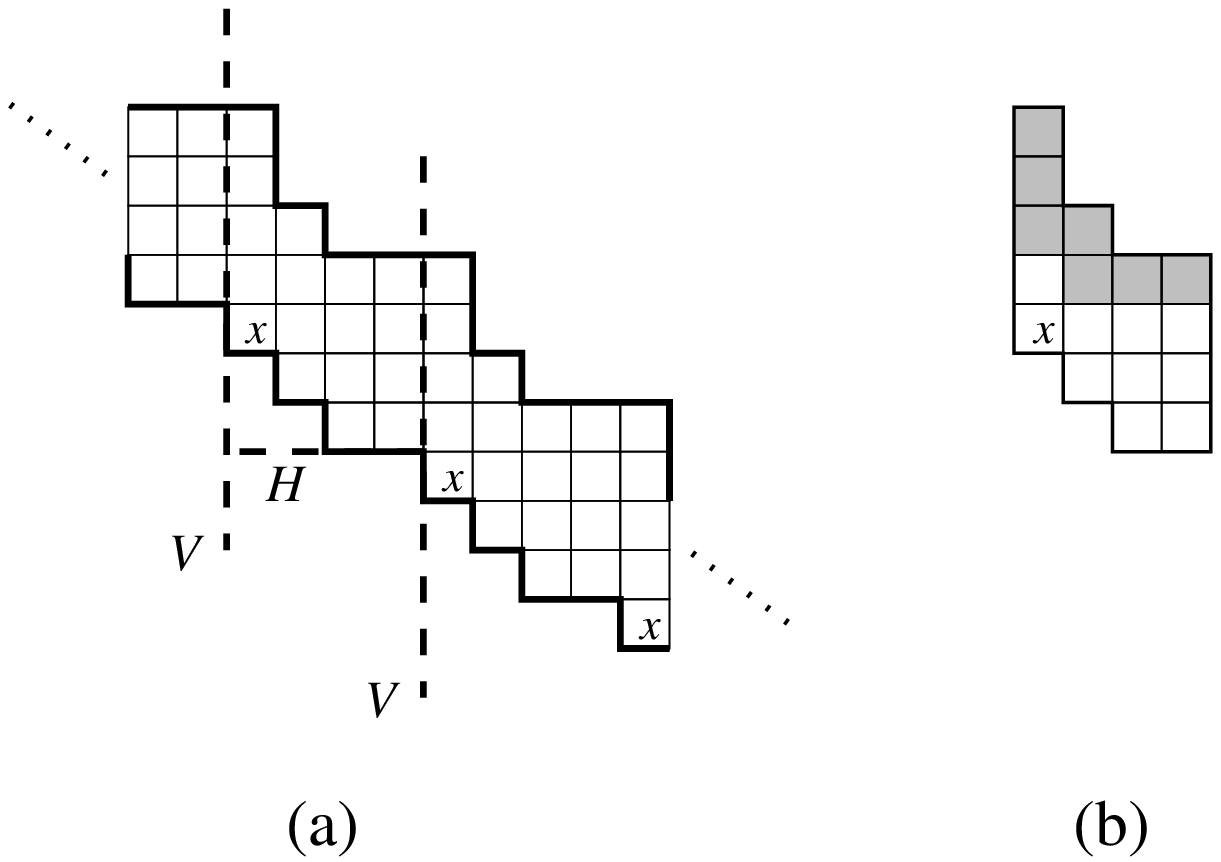} 
\caption{Describing $\C$ as $\ldm$}
\label{fig:findingldm} 
\end{figure}

Suppose we are given any cylindric skew shape $\C$.  The process for finding $\la$,
$d$ and $\mu$ is best understood from a figure, and 
we will use the cylindric skew shape shown in Figure \ref{fig:findingldm}(a) as a 
running example.  The boxes labelled $x$ are identified, so that $k=3$ and 
$n-k=4$ in this example.  First, we must choose a set of representatives for the 
elements of $\C$.  A convenient way to do this is to take the elements
between two adjacent representatives of a vertical line $V$.
Now draw a horizontal line segment $H$ running below each of 
our representatives of $\C$.  In Figure \ref{fig:findingldm}(a),
we regard the intersection of $V$ and the left end of $H$ as our origin.  
The partition $\mu$ is now the
partition whose Young diagram is outlined by $H$, $V$ and the lower boundary of
$\C$.  In our example, $\mu = (2,1)$.  

Next, consider just our set of representatives for the elements of $\C$ 
as in Figure \ref{fig:findingldm}(b).  Define a partition $\La$ by supposing 
the resulting skew shape is $\La/\mu$.  Therefore, in our example, 
$\La = (4,4,4,4,2,1,1)$.  If $\La \inpkn$ 
then set $d=0$, $\la=\La$ and 
we are done.  Otherwise, let $\Lam$ denote the unique
partition $\nu$ that makes $\La/\nu$ an $n$-ribbon with
$n-k$ non-empty columns.  In other words, $\Lam$ is 
obtained by removing an $n$-ribbon along the top of $\La$, starting
in $\La$'s leftmost column and ending in $\La$'s rightmost column.
It is not difficult to see that such a ribbon always has $k+1$ non-empty 
rows. 
In our example, we remove the shaded boxes in Figure \ref{fig:findingldm}(b)
and 
$\Lam = (4,4,4,1)$.  We can see that $\Lam$
is well-defined by referring back to Figure \ref{fig:findingldm}(a).  Effectively
what we are doing is removing the cylindric ribbon that runs all the way along
the top of $\C$.  We see that this cylindric ribbon must have $n$ elements.

Now if $\Lam \inpkn$, then we set $d=1$ and $\la = \Lam$.
Otherwise, obtain $\addribbon{\La}{-2}$ from $\Lam$ in the same 
way that $\Lam$ was obtained from $\La$: remove an $n$-ribbon 
from the top of $\Lam$, 
starting in the leftmost column and ending in the rightmost
column.  
Repeating this procedure, we can construct $\addribbon{\La}{-e}$, stopping as
soon as $\addribbon{\La}{-e} \inpkn$.  
We then set $d=e$ and $\la=\addribbon{\La}{-e}$.
In our example, we see that $\addribbon{\La}{-2} = (3,3) \inpkn$ and so $d=2$,  
$\la = (3,3)$ and $\ldm = (3,3)/2/(2,1)$.

\begin{remark}\label{rem:ldm}  There are several things to note about $\ldm$:
\begin{enumerate}
\item[(i)] For a given $\C$, $\ldm$ is clearly not unique and depends on our choice of 
origin.
\item[(ii)]  $\mu$ is not necessarily contained in $\la$.  For example, moving our 
origin 1 square down and 1 square to the left, the reader is encouraged to verify that
$\ldm = (3,3)/3/(4,3,2,1)$.  This is an example of the following more general statement.
Suppose we have a cylindric shape $\ldm$ with $d\geq1$ and $\mu$ is a partition for which
$\addribbon{\mu}{-1}$ exists.  Then $\ldm$ is the same cylindric shape as
$\la/(d-1)/\addribbon{\mu}{-1}$.
\item[(iii)]
We always have $\la \inpkn$ and it is always possible to choose our origin so
that $\mu \inpkn$.  
\item[(iv)] If $\tau = \mu$, $\La$ or $\addribbon{\La}{-i}$ for $1 \leq i \leq e$, then
$\tau$ satisfies
\[
\tau'_1 \geq \tau'_2 \geq \cdots \geq \tau'_{n-k} \geq \tau'_1 - k. 
\]
\item[(v)] We could alternatively have defined $\la$ by saying it is the
\emph{$n$-core} of $\La$, where the $n$-core is defined in the following
manner.  Given a partition $\tau$, successively remove $n$-ribbons from 
$\tau$ so that after each ribbon removal, the resulting shape is a partition. 
Stop when no more $n$-ribbons can be removed.  It is a well-known fact
(see, for example, \cite[I.1, Example 8]{Mac95}) that the resulting partition $\la$
is independent of the choice of ribbons removed, and $\la$ is 
said to be the \emph{$n$-core} of $\tau$.  
\item[(vi)]
Our notation $\ldm$ is equivalent to that in \cite{Pos04pr}, but our 
explanation of it is very different.  We choose this description in terms of removal of
ribbons because it will be useful in later sections.  
\end{enumerate}
\end{remark}

For any formal power series $f$ in the variables $x = (x_1, x_2, \ldots)$, we will
write $f(x_1,  \ldots, x_k)$ to denote the specialization $f(x_1, x_2, \ldots, x_k, 0, 0, \ldots)$.
We are finally ready to state \cite[Theorem 6.3]{Pos04pr}.

\begin{theorem}\label{thm:postnikov}
For any two partitions $\la, \mu \inpkn$ and a non-negative integer $d$, we have
\begin{equation}\label{equ:postnikov}
s_{\ldm}(x_1, \ldots, x_k) = \sum_{\nu \inpkn} \gwldmn s_\nu (x_1, \ldots, x_k) .
\end{equation}
\end{theorem}

Since we are restricting to $k$ variables, the left-hand side is a sum over semistandard
cylindric tableaux $T$ that map $\ldm$ to the set $[k]$.  Since $T$ must increase in the columns
of $\ldm$, this implies that $s_{\ldm}(x_1, x_2, \ldots, x_k)$ is non-zero only if all
the columns of $\ldm$ contain at most $k$ elements.  One can check that this is equivalent to 
all the rows of $\ldm$ containing at most $n-k$ elements.  In this case, we follow Postnikov
in saying that $\ldm$ is a \emph{toric} shape.  While we take this opportunity to note that 
toric shapes are the shapes that are most
relevant to the Gromov-Witten invariants, we will continue to work with general cylindric
skew shapes.  

While we will be mostly interested in the case of infinitely many variables
\linebreak
$x = (x_1, x_2, \ldots)$, we make a few quick remarks about both
$s_{\ldm}(x)$ and 
\linebreak
$s_{\ldm}(x_1, \ldots, x_k)$.
First, since all the entries in any column of a semistandard cylindric tableau are distinct,
the monomial $x_1^{a_1} x_2^{a_2} \cdots$ appears with coefficient $0$ in 
$s_{\ldm}(x)$ if $a_i > n-k$ for some $i$.  This gives the useful fact that
\begin{equation}\label{equ:partsizes}
s_{\ldm}(x) = \sum_{\nu} c_\nu s_\nu(x) = \sum_{\nu:\nu_1 \leq n-k} c_\nu s_\nu(x).
\end{equation}
From this, we conclude
\[
s_{\ldm}(x_1, \ldots, x_k) = \sum_{\nu: l(\nu) \leq k} c_\nu s_\nu(x_1, \ldots, x_k) 
= \sum_{\nu \inpkn} c_\nu s_\nu(x_1, \ldots, x_k),
\]
explaining why the sum in \eqref{equ:postnikov} is only over $\nu \inpkn$.
Finally, we note that 
$s_{\ldm}(x_1, \ldots, x_k)$ is essentially obtained from $s_{\ldm}(x)$ by 
removing all those terms involving $s_\nu$ with $l(\nu) > k$.  In fact, in the
sections that follow, we will be focusing
most of our attention on these terms $s_\nu$ with $l(\nu) > k$.  

Since we know from the geometric definition of Gromov-Witten invariants that 
$\gwldmn \geq 0$, we conclude that $s_{\ldm}(x_1, \ldots, x_k)$ is Schur-positive.  
On the other hand, we observe that $s_{\ldm}(x)$ may not be
Schur-positive.  For example, the cylindric skew shape $\C$ from Example
\ref{exa:bowtie} has 
\[
s_{\C} = m_{22} + 2m_{211} + 4m_{1111} = s_{22} + s_{211} - s_{1111}.
\]
In the next section, we answer the following question:

\begin{question}
For what cylindric skew shapes $\C$ is
$s_{\C}(x)$ Schur-positive?  
\end{question}


\section{Schur-positivity}\label{sec:schurpositivity}

Recall that, unless otherwise stated, all of our cylindric skew shapes $\C$ 
will be subposets of the cylinder $\cnk$.  
We saw in Example \ref{exa:skewshapecylindric} that the skew shape $\la/\mu$ can 
be regarded as a cylindric skew shape $\C$ when $\la/\mu$ fits inside a box of height
$k$ and width $n-k$.  In this case, we then know that $\cssf{\C}$ is Schur-positive.  
The following  theorem, which is the main result of this section, states that these
are the only Schur-positive cylindric skew Schur functions.  Recall that every cylindric 
skew shape can be viewed as an oriented poset, and it will be
convenient to use this viewpoint for the first half of this section.  We will say that 
two cylindric skew shapes are isomorphic if their corresponding oriented posets are
isomorphic.

\begin{theorem}\label{thm:spositivity}
Let $\C$ be a cylindric skew shape.  Then $\cssf{\C}(x)$ is Schur-positive if and only if 
$\C$ is isomorphic to a skew shape.  
\end{theorem}

In other words, $\cssf{\C}$ is never Schur-positive except in the trivial case of $\C$ being
a skew shape.  As we will see in Theorem \ref{thm:fpositivity}, the same result applies
with ``Schur-positive'' replaced by ``$F$-positive.''

Before proving Theorem \ref{thm:spositivity}, we consider the Schur expansion of 
cylindric ribbons.  While this example is interesting itself, 
it will also play a key role in the proof of Theorem \ref{thm:spositivity}.  
We will identify the cylindric ribbon
$\C$ with
its corresponding oriented poset $(P,\Or)$, enabling us to talk about weak and strict 
``edges'' of $\C$.  In particular, $\C$ must have
$n$ elements, $k$ strict edges, and $n-k$ weak edges.  
We begin with a special class of cylindric ribbons.

\begin{example}
A cylindric ribbon is said to be a \emph{cylindric hook} if it has a unique minimal element
(when viewed as an oriented poset).  See Figure \ref{fig:cylindrichook} for an example.
\begin{figure}
\center 
\epsfxsize=75mm
\epsfbox{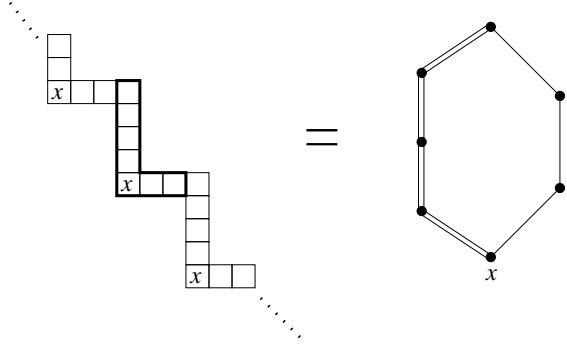} 
\caption{The cylindric hook $\chook{3}{4}$}
\label{fig:cylindrichook} 
\end{figure}
We see that, unlike hooks in the classical case, 
cylindric hooks have just one maximal element.  Also note that $\cnk$
has just one cylindric hook as a subposet, up to isomorphism.  We denote this
cylindric hook by $\hnk$.  Cylindric hooks are the simplest example of a cylindric skew
shape that is not toric.  It follows that $\cssf{\hnk}(x_1, \ldots, x_k)=0$.
This is also evident in the following result which shows that the Schur expansion of 
$\cssf{\hnk}(x)$ is a nice alternating sum of Schur functions of hooks.
\end{example}

\begin{lemma}\label{lem:cylindrichook}

With all functions in the variables $x=(x_1, x_2, \ldots)$, we have
\[
\cssf{\hnk} = s_{(n-k, 1^k)} - s_{(n-k-1, 1^{k+1})} + \cdots +(-1)^{n-k-2}s_{(2, 1^{n-2})} 
+ (-1)^{n-k-1}s_{(1^n)}.
\]
\end{lemma}

We will be ready to prove this lemma as soon as we have introduced a basic tool
that will be important for dealing with cylindric ribbons.  
Suppose $(P,\Or)$ is 
an oriented poset with two incomparable elements $y$ and $z$.  In a 
\pop\ $f$,
either $f(y) \leq f(z)$ or $f(z) < f(y)$.  Let $P(y \leq z)$ denote the oriented 
poset obtained from $(P,\Or)$
by inserting a weak edge from $y$ up to $z$, and let $P(z<y)$ denote the oriented
poset obtained from $(P,\Or)$ by inserting a strict edge from $z$ up to $y$.  
Finally, let us write $P(y \parallel z)$ for the oriented
poset $(P,\Or)$.
We therefore have
\begin{equation}\label{equ:deletionreversal}
K_{P(y \parallel z)}(x) = K_{P(y \leq z)}(x) + K_{P(z < y)}(x).
\end{equation}
For the sake of legibility, 
we will sometimes write or draw $(P,\Or)$ in place of $\kpox$
so that \eqref{equ:deletionreversal} becomes
\[
P(y \parallel z) = P(y \leq z) + P(z < y).
\]
The pair of equations below then follow, and we will refer to them as the 
\linebreak
``deletion-minus-reversal rule'':
\begin{eqnarray}
P(y \leq z) & = & P(y \parallel z) - P(z < y) ,\nonumber \\
P(z < y) & = & P(y \parallel z) - P(y \leq z) . \label{equ:delminusrev}
\end{eqnarray}
\begin{figure}
\center 
\epsfxsize=75mm
\epsfbox{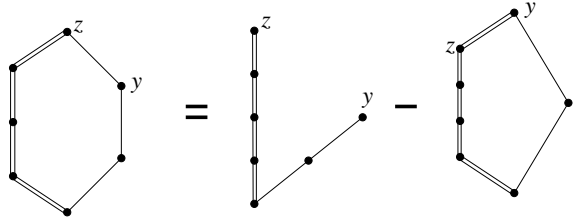} 
\caption{A demonstration of the deletion-minus-reversal rule}
\label{fig:deletionminusreversal} 
\end{figure}
To see this rule in action, see Figure \ref{fig:deletionminusreversal}.  We pick the weak edge
$(y,z)$ in the leftmost poset $P(y \leq z)$ as shown.  Deleting this edge, we get the middle
oriented poset $P(y \parallel z)$.  Reversing the edge and making it strict 
gives the oriented poset $P(z < y)$ on the right.  The deletion-minus-reversal rule gives
an equation among the generating functions, as represented in the figure.  In this particular
case, we get $\chook{3}{4}=(3,1^4) - \chook{2}{5}$.

\begin{proof}[Proof of Lemma \ref{lem:cylindrichook}]
With $n$ fixed, we prove the result by induction on $n-k$, the number of weak edges of $\hnk$.

$\chook{1}{n-1}$ consists of a chain of $n$ elements with $n-1$ strict edges.  (The weak edge
that goes from the bottom element to the top element is redundant and hence is discarded.
Compare this with Note \ref{not:uvequal1}.)
Therefore, $\cssf{\chook{1}{n-1}} = s_{(1^n)}$, as required.  

By the deletion-minus-reversal rule applied to the uppermost weak edge of
\linebreak
$\hnk$, we get
that 
\begin{equation}\label{equ:hooks}
\hnk = (n-k, 1^k) - \chook{n-k-1}{k+1},
\end{equation}
and the result follows.  
\end{proof}

\begin{remark}
We saw in the above proof that $\chook{1}{n-1}$, which is a cylindric ribbon, has a Schur-positive 
generating function $\cssf{\chook{1}{n-1}}=s_{(1^n)}$.  This is, however, not a contradiction with 
Theorem \ref{thm:spositivity}, since $\chook{1}{n-1}$ is isomorphic to the skew shape $(1^n)$.
\end{remark}

We are now ready to discuss the Schur expansions of general cylindric ribbons.  

\begin{proposition}\label{pro:ribbon}
Let $\C$ by a cylindric ribbon which is a subposet of the cylinder $\cnk$.
Then
\[
\cssf{\C}(x) = \left( \sum_{\nu \inpkn} c_\nu s_\nu(x) \right) + \cssf{\hnk}(x),
\]
with $c_\nu$ a non-negative integer for all $\nu \inpkn$.
\end{proposition}

\begin{proof}
From \eqref{equ:partsizes}, we know that 
\[
s_{\C}(x) = \sum_{\nu \inpkn} c_\nu s_\nu(x) + 
\sum_{\genfrac{}{}{0pt}{}{\nu:\nu_1 \leq n-k}{ l(\nu) > k}} c_\nu s_\nu(x).
\]
Restricting to $k$ variables eliminates the second sum, and applying Theorem 
\ref{thm:postnikov} then gives that $c_\nu$ is a non-negative
integer for $\nu \inpkn$.

It remains to show that the terms $s_\nu$ in the Schur expansion of $s_{\C}$ that
have $l(\nu) > k$ correspond to the Schur expansion of $\cssf{\hnk}$.  
With $n$ considered fixed, we proceed by induction on $k$, the number of strict 
edges of $\C$.  Like in the previous proof, the base case is somewhat anomalous.
If $k=1$, then $\C$ is already the cylindric hook $\chook{n-1}{1}$ and we are done.
While $\chook{n-1}{1}$ cannot be expressed as an oriented poset, this does not
affect the rest of the proof.  (Again, compare with Note \ref{not:uvequal1}.)

For $k > 1$, we pick a strict edge of $\C$ and apply the deletion-minus-reversal rule to it.
We get that
\[
\cssf{\C}(x) = s_{\la/\mu}(x) - \cssf{\D}(x),
\]
where $\la/\mu$ is a (classical) ribbon with $k-1$ strict edges, and $\D$ is a cylindric
ribbon with $k-1$ strict edges and $n-k+1$ weak edges.  
Applying the induction hypothesis, we have
\begin{eqnarray*}
\cssf{\C}(x) & = & s_{\la/\mu}(x) - \left( \sum_{\nu \subseteq (k-1) \times (n-k+1)} d_\nu s_\nu(x) \right) 
- \cssf{\chook{n-k+1}{k-1}}(x) \\
& = & s_{\la/\mu}(x) - \left( \sum_{\nu \subseteq (k-1) \times (n-k+1)} d_\nu s_\nu(x) \right) 
- s_{(n-k+1, 1^{k-1})}(x) + \cssf{\hnk}(x) 
\end{eqnarray*}
with the second equality coming from \eqref{equ:hooks}.  Since $\la/\mu$ has $k-1$ strict edges,
it has $k$ rows.  Therefore, any term $s_\nu$ in its Schur expansion has at most $k$ rows.  
We conclude that the terms $s_\nu$ in the Schur expansion of $s_{\C}$ that
have $l(\nu) > k$ are exactly the terms from the expansion of $\cssf{\hnk}(x)$, as required.
\end{proof}

\begin{remark}\label{rem:ribbonspositive}
Given that Schur functions are those skew Schur functions that come from skew shapes
with a unique minimal element, let us say that \emph{cylindric Schur functions} are
those cylindric skew Schur functions that come from cylindric skew shapes with a 
unique minimal element.  Now
let $\C$ be a cylindric ribbon which is a subposet of $\cnk$.  Theorem \ref{thm:spositivity}
tells us that $s_{\C}(x)$ is not Schur-positive.  However, Proposition \ref{pro:ribbon}
says that $s_{\C}(x)$ can be expanded as a positive integer linear combination of 
cylindric Schur functions.  Each of these cylindric Schur functions comes from a cylindric
skew shape that is also a subposet of $\cnk$.  In this case, let us say that
$s_{\C}$ is \emph{cylindric Schur-positive}.  Cylindric Schur-positivity will be the subject of
Section \ref{sec:cylschurpos}.
\end{remark}

For our proof of Theorem \ref{thm:spositivity}, it will be helpful to follow 
\cite{MalThesis,MaRe95}
in defining a coproduct for the ring $\qsym$ of quasisymmetric
functions. 
Let $\mathbb{P}'$ denote the set $\{1', 2', \ldots\}$ with the total order $1' < 2' < \cdots$.
Totally order the disjoint union $\mathbb{P} \cup \mathbb{P'}$ by setting $i < j'$ 
for all $i \in \mathbb{P}, j' \in \mathbb{P'}$.  Given a labelled poset $(P,\omega)$, 
suppose we consider \pwp s $\si$ that are maps from $P$ to $\mathbb{P} \cup \mathbb{P}'$,
rather than from $P$ to $\mathbb{P}$.  Letting $y$ denote the set of variables 
$y=(y_1, y_2, \ldots)$, we can then set
\[
K_{P,\omega}(x,y) = \sum_{\si \in \apw} x_1^{\#\si^{-1}(1)}x_2^{\#\si^{-1}(2)}\cdots
y_1^{\#\si^{-1}(1')}y_2^{\#\si^{-1}(2')}\cdots .
\]

Suppose $\al=(\al_1, \ldots, \al_k)$ is a composition of $N$.  
It is not difficult to find a labelled poset $(P,\omega)$ such that 
$\kpwx = F_\al(x)$.  Indeed,
we let $P$ be a chain of elements $p_1 < p_2 < \cdots < p_N$.
Letting $A_i = \sum_{j=1}^{i} \al_j$, we choose $\omega$ 
so that the edge from 
$p_{A_i}$ to $p_{A_i+1}$ is strict for $i=1,\ldots,k-1$, while
all other edges are weak.  
For compositions $\delta=(\delta_1, \ldots, \delta_d)$ and 
$\epsilon=(\epsilon_1, \ldots, \epsilon_e)$ we let 
$\delta\epsilon$ denote the concatenation $(\delta_1, \ldots, \delta_d, 
\epsilon_1, \ldots, \epsilon_e)$, while $\delta\circ\epsilon$ will 
denote the overlap $(\delta_1, \ldots, \delta_{d-1}, \delta_d + \epsilon_1,\epsilon_2, \ldots, \epsilon_e)$.
We can check that
\begin{equation}\label{equ:coproductf}
F_\al(x,y) = \sum_{\delta, \epsilon:\delta\epsilon=\al} F_\delta(x) F_\epsilon(y)
+ \sum_{\delta, \epsilon: \delta\circ\epsilon=\al} F_\delta(x) F_\epsilon(y).
\end{equation}

Since the set $\{F_\al\}$ forms a basis for $\qsym$, it follows that for every quasisymmetric 
function $G(x)$, we can express $G(x,y)$ as a finite sum
\[
G(x,y) = \sum_i G_i(x) G'_i(y),
\]
where $G_i$ and $G'_i$ are themselves quasisymmetric.  This allows us to 
define the \emph{outer coproduct} $\copr: \qsym \to \qsym \otimes \qsym$ by
\[
\copr(G) = \sum_i G_i \otimes G'_i .
\].

If $(P,\Or)$ is an oriented poset and $Q$ is a convex subposet of $P$, we denote
the designation $\Or$ restricted to the edges of $Q$ by $\Or |_Q$.
It follows from our definition of $\copr$ that 
\begin{equation}\label{equ:coproductkpo}
\copr(K_{P,\Or}) = \sum K_{I, \Or |_{I}} \otimes K_{J, \Or |_{J}}
\end{equation}
where the sum is over all disjoint unions $I \cup J$ such that $I$ is an order ideal
of $P$ and $J$ is an order filter (i.e. dual order ideal) of $P$.  
In particular, 
\begin{equation}\label{equ:coproductschur}
\copr(s_\la) = \sum_{\mu \subseteq \la} s_{\mu} \otimes s_{\la/\mu}.
\end{equation}
Thus the outer coproduct for $\qsym$ is just an extension of the outer coproduct
for symmetric functions of \cite{Gei77,Thi91,Zel81}.  
As one might expect, we say that a coproduct is Schur-positive 
(resp. $F$-positive) if
it can be written as linear combination of terms of the form 
$s_{\mu} \otimes s_{\nu}$ (resp. $F_{\al} \otimes F_{\be}$) with all coefficients positive.  

\begin{proof}[Proof of Theorem \ref{thm:spositivity}]
Suppose $\C$ is a cylindric skew shape that is a subposet of the cylinder $\cnk$.
If $\C$ is isomorphic to a skew shape, then we know by the Littlewood-Richardson
rule that $\cssf{\C}(x)$ is Schur-positive.  Now suppose that $\C$ is a cylindric skew
shape that is not isomorphic to a skew shape. 
We note that if $n-k=1$, then $\C$ is isomorphic to a skew shape, so we assume
that $n-k \geq 2$.
We see from \eqref{equ:coproductschur} that the coproduct of 
a Schur-positive function is Schur-positive.  
Our approach will be to show that $\copr(\cssf{\C})$ is not Schur-positive and,
therefore, it will follow that 
$\cssf{\C}(x)$ is not Schur-positive.

Since $\C$ is not isomorphic to a skew shape, $\C$ contains a cylindric
ribbon.  Let $\toprib$ denote the cylindric ribbon with $n$ elements
that runs all the way along the
top of $\C$ and let $\addribbon{\C}{-1}$ denote the
cylindric skew shape that remains after we remove $\toprib$ from $\C$.  
Clearly, viewing $\C$ as an oriented poset, the elements of 
$\addribbon{C}{-1}$ correspond to an order ideal 
of $\C$ and the elements of $\toprib$ correspond
to an order filter of $\C$.  
Choose any partition $\la$ such that
$s_\la$ appears with non-zero coefficient $m$ in $\cssf{\addribbon{C}{-1}}$.
By Proposition \ref{pro:ribbon} and Lemma \ref{lem:cylindrichook}, 
we know that 
\[
[s_{(1^n)}]\cssf{\toprib} = (-1)^{n-k-1} = -[s_{(2,1^{n-2})}]\cssf{\toprib}.
\]
We will now show that
\begin{equation}\label{equ:coproductnegative}
[s_\la \otimes s_{(1^n)}]\copr(\cssf{\C}) = (-1)^{n-k-1}m = -[s_\la \otimes s_{(2,1^{n-2})}]\copr(\cssf{\C}),
\end{equation}
implying that $\copr(\cssf{\C})$ cannot be Schur-positive.  Indeed, suppose $J \neq \toprib$
is an order filter
of $\C$ with $n$ elements.  The only order filter of $\C$ with $n$ elements
that contains a cylindric ribbon is $\toprib$.  Therefore, $J$ does not contain a cylindric ribbon 
and so is isomorphic to a skew shape.  However, any skew shape that is a subposet of 
$\C$ has at most $k$ rows.  Since $k \leq n-2$, we conclude that 
$[s_{(1^n)}]\cssf{J} = [s_{(2,1^{n-2})}]\cssf{J}=0$.  Applying \eqref{equ:coproductkpo},
we now deduce \eqref{equ:coproductnegative}.
\end{proof}

We should justify our earlier assertion that the following result is also true:

\begin{theorem}\label{thm:fpositivity}
Let $\C$ be a cylindric skew shape.  Then $\cssf{\C}(x)$ is $F$-positive if and only if 
$\C$ is isomorphic to a skew shape.  
\end{theorem}

\begin{proof}
Our proof is largely the same as the proof of Theorem \ref{thm:spositivity}.  As is known 
(see, e.g., \cite[Theorem 7.19.7]{ec2}), Schur functions have non-negative coefficients when
expressed in the basis of fundamental quasisymmetric functions $F_\al$.  More specifically,
we define a \emph{standard Young tableau} (SYT) $T$ of shape $\la$ to be a filling
of the Young diagram of $\la$ with distinct entries from the set $\{1, 2, \ldots, |\la |\}$ that
increases in the rows and up the columns (using French notation).  The descent set
of $T$ is defined to be those numbers $i \in \{1,\ldots, |\la |-1\}$ such that 
$i+1$ is in a strictly higher row of $T$ than $i$.  The composition $\mathit{co}(T)$ is
then given by $(i_1, i_2-i_1, \ldots, i_k-i_{k-1}, |\la|-i_k)$, where $\{i_1, \ldots, i_k\}$ is the
descent set of $T$.  We then have
\begin{equation}\label{equ:stof}
s_\la = \sum_{T} F_{\mathit{co}(T)},
\end{equation}
where the sum is over all SYT $T$ of shape $\la$.  

In particular, we see that $[F_{(1^n)}]s_\la = 0$ unless $\la = (1^n)$, and
$[F_{(2,1^{n-2})}]s_\la = 0$ unless $\la=(2,1^{n-2})$.
From Proposition \ref{pro:ribbon} and Lemma \ref{lem:cylindrichook}, it follows that
\[
[F_{(1^n)}]\cssf{\toprib} = (-1)^{n-k-1} = -[F_{(2,1^{n-2})}]\cssf{\toprib},
\]
where $\toprib$ is a cylindric ribbon.
We now mimic the proof of Theorem \ref{thm:spositivity} to show that $\copr(\cssf{\C})$ is not
$F$-positive.  However, by \eqref{equ:coproductf}, the coproduct of an $F$-positive function
is $F$-positive, so the result follows.
\end{proof}

We finish our discussion of $F$-positivity by addressing two interesting issues.

\begin{remark}
Schur-positive functions are $F$-positive by \eqref{equ:stof}, but the converse is
not true.  For example, 
\[
F_{31} + F_{13} + F_{211} + F_{112} = s_{31} + s_{211} - s_{22}.
\]
Therefore, Theorem \ref{thm:fpositivity} is seemingly stronger than Theorem \ref{thm:spositivity}.
However, we chose to prove Theorem \ref{thm:spositivity} separately for two reasons.  The first is that
our main subject is cylindric skew Schur functions and Schur-positivity.  The second
reason is that we have been unable to find a symmetric function of the form 
$\kpox$ that is 
\linebreak
$F$-positive but not Schur-positive.  No such examples exist for 
$|P| \leq 7$.  Determining whether or not an example exists might be an interesting
problem.  More generally, we can ask what quasisymmetric functions can be expressed as 
$\kpox$, or even just as $\kpwx$.  Restricting to symmetric functions, we can also ask how
to easily tell when a positive linear combination of Schur functions is equal to a 
skew Schur function.
\end{remark}

\begin{remark}
One might wonder if Theorem \ref{thm:fpositivity} can be extended to functions that aren't
symmetric.  Specifically, one might ask if the following statement is true:

\textit{Let $(P,\Or)$ be an oriented poset.  Then $\kpox$ is $F$-positive if and only if
$(P, \Or)$ is a labelled poset.}

This statement is false, as shown by the example $(P,\Or)$
in Figure \ref{fig:fpositivitycounterex}.
\begin{figure}
\center 
\epsfxsize=25mm
\epsfbox{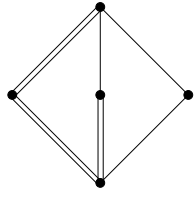} 
\caption{An oriented poset $(P,\Or)$ that is not a labelled poset, but has $\kpox$ $F$-positive}
\label{fig:fpositivitycounterex} 
\end{figure}
It has a cycle, but $\kpox = F_{131} + F_{113} + F_{221} + F_{212} + 2F_{122}$.

This further suggests that, among oriented posets, cylindric skew shapes are noteworthy.
\end{remark}

Let $\C$ be a cylindric skew shape that is not isomorphic to a skew shape.
We know from Theorem \ref{thm:postnikov} that $\cssf{\C}$ in $k$ variables is Schur-positive.
On the other hand, by Theorem \ref{thm:spositivity}, $\cssf{\C}$ in an infinite number of variables 
is not Schur-positive.  We conclude this section with a discussion of the minimum number of 
variables in which $\cssf{\C}$ fails to be Schur-positive.  

As before, let $\toprib$ denote the cylindric ribbon with $n$ elements
that runs all the way along the
top of $\C$ and let $\addribbon{\C}{-1}$ denote the
cylindric skew shape that remains after we remove $\toprib$ from $\C$.  

If $\C$ is a cylindric ribbon, we deduce from Proposition \ref{pro:ribbon} and
Lemma \ref{lem:cylindrichook}  
that $\cssf{\C}$ remains Schur-positive in $k+1$ 
variables but always fails to be 
Schur-positive in $k+2$ variables.  
By looking at coproducts, 
we can use this fact to say something about general cylindric skew shapes. 

\begin{proposition}
Let $\C$ be a cylindric skew shape that is not isomorphic to a skew shape and that is a 
subposet of $\cnk$.
If $m$ denotes the maximum number of elements in a column of $\addribbon{\C}{-1}$, 
then $\cssf{\C}$ is not Schur-positive in $m+k+2$ variables.  
\end{proposition}

\begin{proof}
We begin by finding a partition $\tau$ such that $s_\tau(x)$ appears with positive coefficient 
in the Schur expansion of 
$\cssf{\addribbon{\C}{-1}}(x)$.  We can form a semistandard cylindric tableau $T$ of shape 
$\addribbon{\C}{-1}$ by mapping the $i$th lowest element
of each column to $i$, for all $i$.  
Set $\tau$ to be the \emph{content} of $T$, i.e.,
$\tau = (\#T^{-1}(1), \#T^{-1}(2), \ldots)$.
Notice that $T$ is the only semistandard cylindric tableau of shape 
$\addribbon{\C}{-1}$ and content $\tau$. 
Therefore, when we expand $\cssf{\addribbon{\C}{-1}}$ in terms of the monomial symmetric
functions, $m_\tau$ appears with coefficient $+1$.
Furthermore, we see that $\tau$ is a maximal
possible content of a semistandard cylindric tableau of shape $\addribbon{\C}{-1}$
in \emph{dominance order}.  (This means that if $\si$ is some other possible content, 
then $\sum_{j=1}^{i} \si_j \leq \sum_{j=1}^{i} \tau_j$ for all $i \geq 1$.)
It follows that $s_\tau(x)$ appears with coefficient $+1$ in the Schur expansion of 
$\cssf{\addribbon{\C}{-1}}(x)$.  (If this is not clear, see \cite[Proposition 7.10.5]{ec2}.)

We know that $s_{(n-k-1,1^{k+1})}(x)$ appears with coefficient $-1$ in $\cssf{\toprib}(x)$.
Looking at $\copr(\cssf{\C})$, we now see that $s_\tau \otimes s_{(n-k-1, 1^{k+1})}$
appears with coefficient $-1$.  Comparing this with \eqref{equ:coproductschur},
we see there exists a partition $\la$ such that:
\begin{enumerate}
\item[(i)] $\tau \subseteq \la$, and
\item[(ii)] $s_{(n-k-1, 1^{k+1})}$ appears with positive coefficient in the Schur expansion of
$s_{\la/\tau}$, and
\item[(iii)] $s_\la(x)$ appears with coefficient $-1$ in the Schur expansion of $\cssf{\C}(x)$.
\end{enumerate}
In particular, we know that $l(\la) \leq l(\tau) + l((n-k-1, 1^{k+1})) = m+k+2$.
Therefore, $s_\la(x_1, \ldots, x_{m+k+2}) \neq 0$, and so $s_\la(x_1, \ldots, x_{m+k+2})$ 
appears with coefficient $-1$ in the Schur expansion of 
$\cssf{\C}(x_1, \ldots, x_{m+k+2})$.
\end{proof}

We do not claim, and it is not true, that $m+k+2$ is the best possible value.  In other words, 
it can be the case that $s_C$ is not Schur-positive in some number of variables that is 
less than $m+k+2$.  For toric shapes, it is clear that $m \leq k-1$, and so we get the 
following result.

\begin{corollary}
Let $\C$ be a toric shape that is not isomorphic to a skew shape and that is a 
subposet of $\cnk$.
Then $\cssf{\C}$ is not Schur-positive in $2k+1$ variables.  
\end{corollary}


\section{From cylindric skew shapes to skew shapes}\label{sec:skewexpansion}

So far, essentially the only tool we have for working with cylindric skew Schur functions is 
the deletion-minus-reversal rule of \eqref{equ:delminusrev}.  The subject of this
section is a rule for expressing any cylindric skew Schur function as a signed sum of 
skew Schur functions.  Our rule is based on a result of Gessel and Krattenthaler
from \cite{GeKr97}, with our reformulation modelled on a result from \cite{BCF99}.
We begin with an exposition of these two results, starting with the latter.

By saying that a partition $\tau$ is obtained from $\la$ by adding $d$ $n$-ribbons, or that $\la$ is
obtained from $\tau$ be removing $d$ $n$-ribbons, we mean that
there is a sequence of partitions
\begin{equation}\label{equ:nribbons}
\la = \nu_0 \subseteq \nu_1 \subseteq \cdots \subseteq \nu_d = \tau
\end{equation}
such that $\nu_i/\nu_{i-1}$ is an $n$-ribbon for $i=1,\ldots,d$.  We say that the \emph{width}
of a ribbon is its number of non-empty columns.  If $\tau_1 \leq n-k$, then we define
\[
\vare(\tau/\la) = (-1)^{\sum_{i=1}^{d} (n-k-\mathrm{width}(\nu_i/\nu_{i-1}))}.
\]
It can be shown that $\vare(\tau/\la)$ is independent of the choice of the sequence in
\eqref{equ:nribbons}.

The result of interest from \cite{BCF99} is the following:

\begin{theorem}\label{thm:bcff}
Suppose we have
$\la, \mu, \nu \inpkn$ with $|\mu| + |\nu| = |\la| + dn$ for some $d \geq 0$.
Then the Gromov-Witten invariant $\gwldmn$ can be expressed in terms of 
Littlewood-Richardson coefficients as
\begin{equation}\label{equ:bcff}
\gwldmn = \sum_\tau \vare(\tau/\la) c^{\tau}_{\mu \nu},
\end{equation}
where the sum is over all $\tau$ with $\tau_1 \leq n-k$ that can be obtained from $\la$ by
adding $d$ $n$-ribbons.
\end{theorem}

Formulas for $\gwldmn$ similar to \eqref{equ:bcff} have
appeared in different contexts in \cite{Cum91,GoWe90,Kac90,Wal90}.
Combining Theorems \ref{thm:postnikov} and \ref{thm:bcff}, we get:

\begin{corollary} \label{cor:bcff}
For any cylindric skew shape $\ldm$ with $\la, \mu \inpkn$, we have
\begin{equation}\label{equ:bcff2}
s_{\ldm}(x_1, \ldots, x_k) = \sum_\tau \vare(\tau/\la) s_{\tau/\mu} (x_1, \ldots, x_k),
\end{equation}
where the sum is over all $\tau$ with $\tau_1 \leq n-k$ that can be obtained from $\la$ by
adding $d$ $n$-ribbons.
\end{corollary}

\begin{proof}
Multiply both sides of \eqref{equ:bcff} by $s_\nu (x_1, \ldots, x_k)$, sum over all 
$\nu \inpkn$, and
apply Theorem \ref{thm:postnikov}.
\end{proof}

From our point of view, the obvious disadvantage of Corollary \ref{cor:bcff} is that it 
only gives certain terms in the expansion of $s_{\ldm}(x)$.  For example, for cylindric shapes
that are not toric, both sides of \eqref{equ:bcff2} will be zero.  Gessel and Krattenthaler's setting
does not have this limitation.  To apply their result to get an expression
for $s_{\ldm}(x)$, we first have some work to do.  Their basic result \cite[Proposition 1]{GeKr97} is stated
in terms of lattice paths.  In \cite[\S 9]{GeKr97} , they show how to apply Proposition 1 to 
obtain expressions for Schur functions.  Mimicking their approach, we first obtain an
expression for $s_{\ldm}$ in terms of the elementary symmetric functions.  Recall from
our construction preceding Remark \ref{rem:ldm} that, for a given $\ldm$, $\La$ is the unique partition 
satisfying $\addribbon{\La}{-d} = \la$.  In this case, we also write $\addribbon{\la}{d} = \La$
and we see that $\La$ is obtained from $\la$ by adding $d$ $n$-ribbons, each starting in 
$\la$'s rightmost column (column $n-k$) and ending in column 1. 
We get that
\footnote{
For the benefit of the reader wishing to derive \eqref{equ:gkbasic} from \cite{GeKr97}, we note 
that we took $z=1$, $\textbf{S}=(-n,n)$, $w(e) = x_{i+j}$, $n=m$,
$u_i = (-(\La'_i+n-k-i), \La'_i+n-k-i)$, 
\linebreak
$v_i = (-(\mu'_i+n-k-i), \mu'_i+n-k-i+m)$, $r_i=-k_i$, and we let $m$ tend to infinity.}
\begin{equation}\label{equ:gkbasic}
s_{\ldm}(x) = \sum_{\genfrac{}{}{0pt}{}{r_1 + \cdots +r_{n-k}=0}{r_i \in \mbz}}
\det \left( e_{r_s n + \La'_s - \mu'_t -s + t}(x) \right)_{s,t = 1}^{n-k}.
\end{equation}
As usual, we set $e_0 = 1$ and $e_i = 0$ for $i < 0$.  
The alert reader may notice the possibility of greatly simplifying \eqref{equ:gkbasic} using
the dual Jacobi-Trudi
identity (see 
\linebreak 
\cite[I, (5.5)]{Mac95} or \cite[Corollary 7.16.2]{ec2}):
\[
s_{\tau'/\mu}(x) = \det \left( e_{\tau_s - \mu'_t -s + t}(x)\right)_{s,t=1}^{n-k},
\]
where $\mu \subseteq \tau'$ and $\tau_1 \leq n-k$.
Indeed, given $r=(r_1, \ldots, r_{n-k})$, 
let $\modif{\La'}{r}$ denote the integer sequence $(\La'_1+r_1n, \ldots, \La'_{n-k}+r_{n-k}n)$.
Now $\modif{\La'}{r}$ may not be a partition.  However, we can still define the 
Schur function $s_{\al'/\mu}$ for any sequence $\al=(\al_1, \ldots, \al_{n-k})$ using 
a dual Jacobi-Trudi determinant:
\begin{equation}\label{equ:generaljt}
s_{\al'/\mu} = \det \left( e_{\al_s - \mu'_t -s + t}\right)_{s,t=1}^{n-k}.
\end{equation}
By repeatedly transposing adjacent rows, it may be possible to make the matrix on
the right-hand side of \eqref{equ:generaljt} into the dual Jacobi-Trudi matrix of a 
skew shape $\tau'/\mu$, multiplied by a sign term $(-1)^{\delta(\al)}$.
For example,
\begin{eqnarray*}
s_{(7,-2,4,11)'/(2,1,0,0)} & = &
\left| \begin{array}{cccc}
e_{5} &  e_{7} & e_{9} & e_{10} \\
e_{-5} &  e_{-3} & e_{-1} & e_{0} \\
e_{0} &  e_{2} & e_{4} & e_{5} \\
e_{6} &  e_{8} & e_{10} & e_{11} 
\end{array} \right| = +
\left| \begin{array}{cccc}
e_{6} &  e_{8} & e_{10} & e_{11} \\
e_{5} &  e_{7} & e_{9} & e_{10} \\
e_{0} &  e_{2} & e_{4} & e_{5} \\
e_{-5} &  e_{-3} & e_{-1} & e_{0} 
\end{array} \right|  \\
& = & s_{(8,8,4,0)'/(2,1,0,0)}.
\end{eqnarray*}
To save us having to always think in terms of determinants, we can view this process
another way.
Effectively what we are doing is defining an equivalence relation
$\sim$ on integer sequences by saying that 
\[
(\al_1, \al_2, \ldots, \al_{n-k}) \sim
(\al_1, \ldots, \al_{i-1}, \al_{i+1}-1, \al_i+1, \al_{i+2}, \ldots, \al_{n-k}):
\]
we transpose two adjacent elements
of the sequence, increasing the element moving right by 1 and decreasing the element moving
left by 1.  We see that every equivalence class of a sequence $\al$ contains at most one 
partition $\tau$.  If $\tau'$
contains $\mu$ then we
say that $s_{\al'/\mu} = (-1)^{\delta(\al)} s_{\tau'/\mu}$ whenever $\al \sim \tau$, where
$\delta(\al)$ is the number of adjacent transpositions necessary to make $\al$ into a partition.
If $\tau'$ does not contain $\mu$ or if $\al$ does not have a partition in its equivalence class,
then we set $s_{\al'/\mu}=0$.  One can check that this is consistent with the definition 
\eqref{equ:generaljt} of $s_{\al'/\mu}$ as a determinant.
In our example above, we would have had
\[
\al = (7, -2, 4, 11) \sim (7, 3, -1, 11) \sim (7, 3, 10, 0) \sim (7, 9, 4, 0) \sim (8, 8, 4, 0) = \tau
\]
and $\delta(\al)=4$.

Putting this all together, \eqref{equ:gkbasic} becomes:
\begin{theorem}\label{thm:gk}  \emph{\cite{GeKr97} }
For any cylindric shape $\ldm$ that is a subposet of $\cnk$,
we have
\[
s_{\ldm}(x) = \sum_{\genfrac{}{}{0pt}{}{r_1 + \cdots +r_{n-k}=0}{r_i \in \mbz}}
s_{(\modif{\La'}{r})'/\mu}(x)
\]
where $\La = \addribbon{\la}{d}$.  
\end{theorem}

\begin{example}\label{exa:gk}
Consider $\ldm=(3,3)/2/(2,1)$ as depicted in Figure \ref{fig:findingldm}.  
We see that $n=7$, $n-k=4$, $\La'=(7,5,4,4)$ and $\mu=(2,1,0,0)$.  The
values of $r=(r_1, \ldots, r_{n-k})$ that make $s_{(\modif{\La'}{r})'/\mu} \neq 0$ 
are listed in the first column of Table \ref{tab:gkexample}.
\begin{table} 
\[
\begin{array}{cccc}
r & \modif{\La'}{r} & \tau & \delta(\modif{\La'}{r})  \\
\hline
(0,0,0,0)  & (7,5,4,4) & (7,5,4,4) & 0 \\
(-1, 0, 0, 1) & (0, 5, 4, 11) & (8, 5, 4, 3) & 5 \\
(-1, 0, 1, 0) & (0, 5, 11, 4) & (9, 5, 3, 3) & 4 \\
(-1, 1, 0, 0) & (0, 12, 4, 4) & (11, 3, 3, 3) & 3 \\
(0, -1, 0, 1) & (7,-2,4,11) & (8,8,4,0) & 4 \\
(0, -1, 1, 0) & (7,-2,11,4) & (9, 8, 3, 0) & 3 \\
(1,-1,0,0) & (14,-2,4,4) & (14,3,3,0) & 2 \\
(-1,-1,1,1) & (0,-2,11,11) & (9,9,2,0) & 4 \\
(-1,-1,0,2) & (0,-2,4,18) & (15,3,2,0) & 5 \\
(-1,-1,2,0) & (0,-2,18,4) & (16,2,2,0) & 4
\end{array}
\]
\caption{Applying Theorem \ref{thm:gk} to $\ldm=(3,3)/2/(2,1)$}
\label{tab:gkexample}
\end{table}
We conclude that 
\begin{eqnarray*}
s_{(3,3)/2/(2,1)}(x) & = & s_{(7,5,4,4)'/(2,1)}(x) - s_{(8,5,4,3)'/(2,1)}(x)  +s_{(9, 5, 3, 3)'/(2,1)}(x) \\ & & 
-s_{(11, 3, 3, 3)'/(2,1)}(x)  +s_{(8,8,4,0)'/(2,1)}(x) -s_{(9, 8, 3, 0)'(2,1)}(x) \\ & & 
+s_{(14,3,3,0)'/(2,1)}(x) +s_{(9,9,2,0)'/(2,1)}(x) -s_{(15,3,2,0)'/(2,1)}(x) \\ & & 
+s_{(16,2,2,0)'/(2,1)}(x).
\end{eqnarray*}
\end{example}

Using Theorem \ref{thm:gk}, we can actually show that Corollary \ref{cor:bcff}
extends to the case of infinitely many variables.  This is the main result of
this section.

\begin{theorem}\label{thm:gkribbons}
For any cylindric skew shape $\ldm$ that is a subposet of 
$\cnk$, we have
\begin{equation}\label{equ:gkribbons}
s_{\ldm}(x) = \sum_\tau \vare(\tau/\la) s_{\tau/\mu} (x),
\end{equation}
where the sum is over all $\tau$ with $\tau_1 \leq n-k$ that can be obtained from $\la$ by
adding $d$ $n$-ribbons.
\end{theorem}

\begin{note}
While $\la \inpkn$ by definition, we do not require that $l(\mu) \leq k$, unlike
in Theorem \ref{thm:bcff} and Corollary \ref{cor:bcff}.
\end{note}

\begin{example}
Again, consider $\ldm=(3,3)/2/(2,1)$ as depicted in Figure \ref{fig:findingldm}. 
Figure \ref{fig:gkribbons} shows the set of all possible $\vare(\tau/\la) \tau/\mu$ 
with $\tau_1 \leq n-k$
such that $\tau$ can be obtained from $(3,3)$ by adding two 7-ribbons.
\begin{figure}
\center 
\epsfxsize=60mm
\epsfbox{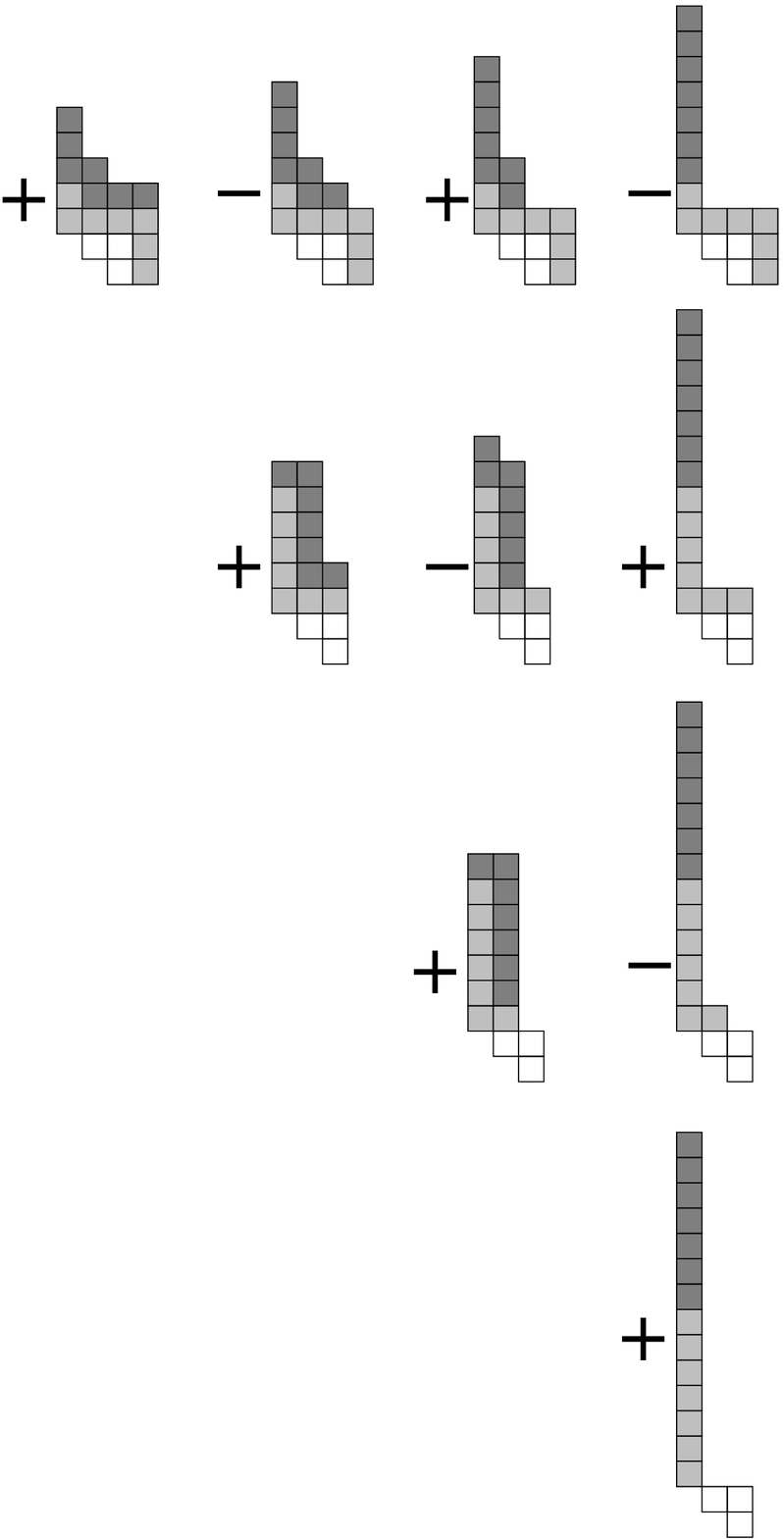} 
\caption{All possible $\vare(\tau/\la) \tau/\mu$ in Theorem \ref{thm:gkribbons}
when $\ldm=(3,3)/2/(2,1)$}
\label{fig:gkribbons} 
\end{figure}
The positioning of the partitions in the figure is supposed to be helpful, as it is determined
by the rightmost column of the added
ribbons.  There can be more than one way to add ribbons to $\la$ and
get a particular $\tau$, but this does not affect our expression for $s_{\ldm}$.

We see that we get the same result as in Example \ref{exa:gk}.  
While the result obtained from Theorem \ref{thm:gk} is more compact to
write, we find the graphical description of $s_{\ldm}$ in Theorem \ref{thm:gkribbons} 
preferable, especially from the point of view of intuition.  We will make
much use of Theorem \ref{thm:gkribbons} in the next section.
\end{example}

\begin{remark}
Because the expression of a cylindric skew shape $\C$ in the form $\ldm$ is not unique, 
Theorem \ref{thm:gkribbons} can be used to give a host of identities among skew Schur
functions.  For example, consider the cylindric skew shape $\C$ shown in Figure 
\ref{fig:identities} with $k=n-k=3$.  
\begin{figure}
\center 
\epsfxsize=35mm
\epsfbox{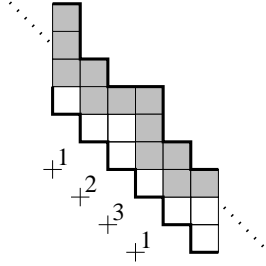} 
\caption{$(3,3,1)/1/(2,1) = (3,2,2)/1/(2,1) =  (1)/2/(2,1)$}
\label{fig:identities} 
\end{figure}
By choosing the origins labelled 1, 2 and 3 respectively, 
we see that $\C$ can be written as $(3,3,1)/1/(2,1)$, $(3,2,2)/1/(2,1)$ or 
$(1)/2/(2,1)$.  It follows that
\begin{eqnarray*}
\cssf{\C}(x) & = & s_{333211/21} - s_{3322111/21} + s_{331111111/21} \\
& = & s_{33331/21} - s_{32221111/21} + s_{322111111/21} \\
& = & s_{33322/21} - s_{3222211/21} + s_{3211111111/21} + s_{2222221/21} - s_{22111111111/21}.
\end{eqnarray*}
\end{remark}

The remainder of this section, which is somewhat technical,
is devoted to working towards and proving 
Theorem \ref{thm:gkribbons}.  
It does not seem that the proof of Theorem \ref{thm:bcff} from \cite{BCF99}
can be easily modified to work in this more general setting.  Instead, our approach will
to show that the statements of Theorems \ref{thm:gk} and \ref{thm:gkribbons}
are equivalent, thereby implying Theorem \ref{thm:gkribbons}.

We begin with some preliminary results about
the $\sim$ equivalence relation.  Rather than working with integer sequences, it 
will be more convenient to work now with \emph{signed} integer sequences.  A signed 
integer sequence is simply an integer sequence with a purely symbolic sign in
front.  By this, we mean that 
 $-(\al_1, \ldots, \al_{n-k})$ is certainly not the same thing as 
$(-\al_1, \ldots, -\al_{n-k})$.  However, we will say that 
$-(-(\al_1, \ldots, \al_{n-k})) = (\al_1, \ldots, \al_{n-k})$.
We extend $\sim$ to the class of signed integer sequences by saying that
\[
(\al_1, \al_2, \ldots, \al_{n-k}) \sim
- (\al_1, \ldots, \al_{i-1}, \al_{i+1}-1, \al_i+1, \al_{i+2}, \ldots, \al_{n-k}),
\]
i.e., the sign changes when we do an adjacent transposition.  
Signed partitions are then defined in the obvious way, and we denote
the set of signed partitions by $\SPar$.  We identify the partition $\la$ with 
the signed partition $+\la$. 

We remark that a signed integer sequence $\pm(\al_1, \ldots, \al_{n-k})$ may not always 
have a signed partition in its equivalence class.  For example, any signed
sequence equivalent
to $-(7,4,0,-2)$ will always have a negative entry in its sequence.  
More interestingly, there is no signed sequence
in the equivalence class of $(7,4,5,0)$ whose sequence is weakly decreasing.
However, we see that any integer sequence 
$\al=(-1)^k(\al_1, \ldots, \al_{n-k})$ is equivalent to a unique signed integer sequence 
$\be=(-1)^{k+\delta(\al)}(\be_1, \be_2, \ldots, \be_{n-k})$ with
\[
\be_1 -1 \geq \be_2 -2 \geq \ldots \geq \be_{n-k} - (n-k).
\]
Here $\delta(\al)$ is the number of adjacent transpositions necessary to
convert $\al$ to $\be$.
We then denote this signed sequence $\be$ by
$\order{(\al_1, \ldots, \al_{n-k})}$ or just $\order{\al_1, \ldots, \al_{n-k}}$.
Finally, if $\al$ is a signed integer sequence, 
we let $\increm{\al}{i}{m}$ denote the signed integer sequence that results
when we increase the $i$th element of $\al$ by $m$, but leave the 
sign of $\al$ unchanged.  Using a similar principle, $\al + (r_1, \ldots, r_{n-k})n$
denotes the signed integer sequence that results when we increase
the $i$th element of $\al$ by $r_i n$ for $i=1, \ldots,n-k$, 
but leave the sign of $\al$ unchanged.

Our first lemma, while only a small portion of the work to come, highlights the
basic connection between ribbons and the $\sim$ equivalence relation.  

\begin{lemma}\label{lem:simtoribbons}
\begin{enumerate}
\item[(a)]
For a partition $\tau$ with $\tau_1 \leq n-k$, suppose we can add an $n$-ribbon
to $\tau$ whose rightmost column is column $i$ to get a new partition $\si$.  
Then $\si$ exists if and only if $\order{\increm{\tau'}{i}{n}}$ is a signed partition, in which
case $\vare(\si/\tau) \si' = (-1)^{n-k-1} \order{\increm{\tau'}{i}{n}}$.
\item[(b)]
For a partition $\tau$ with $\tau_1 \leq n-k$, suppose we can remove an $n$-ribbon
from $\tau$ whose leftmost column is column $i$ to get a new partition $\si$.  
Then $\si$ exists if and only if $\order{\increm{\tau'}{i}{-n}}$ is a signed partition, in which
case 
\linebreak
$\vare(\si/\tau) \si' = (-1)^{n-k-1} \order{\increm{\tau'}{i}{-n}}$.
\end{enumerate}
\end{lemma}

\begin{proof}
We prove (a), with (b) being similar.  
We have that
\begin{eqnarray}
& & (-1)^{n-k-1} \order{\increm{\tau'}{i}{n}}  \nonumber \\
& = &  (-1)^{n-k-1} \order{\tau'_1, \ldots, \tau'_{i-1}, \tau'_i +n, \tau'_{i+1}, \ldots, \tau'_{n-k}} \nonumber \\
& = & (-1)^{n-k-1+i-j} (\tau'_1, \ldots, \tau'_{j-1}, \tau'_i +n-(i-j), \tau'_j+1, \ldots, \tau'_{i-1}+1, \nonumber \\
& & \tau'_{i+1}, \ldots, \tau'_{n-k}),  \label{equ:addedribbonsequence}
\end{eqnarray}
where we take $j$ to be as small as possible subject to the condition that 
\linebreak
$\tau'_i +n-(i-j) \geq \tau'_j+1$.  We observe that \eqref{equ:addedribbonsequence}
gives  exactly the 
column heights of the result $\si$ of adding an $n$-ribbon to $\tau$ whose rightmost
column is $i$ and whose leftmost column is as far left as possible.  Obviously, 
$\si$ is a partition if and only if \eqref{equ:addedribbonsequence} is a signed partition.
(This is the case if and only if 
$\tau'_{j-1} \geq \tau'_i +n-(i-j)$.)
Finally, $\vare(\si/\tau) = (-1)^{n-k-(i-j+1)}$, as required.
\end{proof}

The next lemma encompasses the remaining preliminaries necessary for proving
Theorem \ref{thm:gkribbons}.

\begin{lemma}\label{lem:mult1}
Suppose $\ldm$ with $d \geq 1$ is a cylindric skew shape.  Let  
\linebreak
$r = (r_1, \ldots, r_{n-k})$ be an integer sequence, and let $\La = \addribbon{\la}{d}$.
\begin{enumerate}
\item[(a)] If $\order{\La'+rn} = \order{\La'+tn}$ for some integer 
sequence $t=(t_1, \ldots, t_{n-k})$, then $r=t$.  
\item[(b)] $\order{\La'+rn}$ is a weakly decreasing sequence.
\item[(c)] If $r_1 + \cdots + r_{n-k} = 0$, then $\order{\La'+rn}_1 \geq k+1$.  
\item[(d)] If $r_1 + \cdots + r_{n-k} = 0$ and $\order{\La'+rn}$ is a signed partition, then the signed sequence $\order{\increm{\order{\La'+r n}}{1}{-n}}$ is a signed partition.
\end{enumerate}
\end{lemma}

\begin{proof}
(a) Suppose that $\order{\La'+rn} = \pm(\tau_1, \ldots, \tau_{n-k})$.
Then for some set
\linebreak
$\{i_1, \ldots, i_{n-k}\} = \{1, \ldots, n-k\}$
we have the following congruences
modulo $n$:

\begin{eqnarray*}
\tau_1 & \equiv & \La'_{i_1} + (1-i_1) \\
\tau_2 & \equiv & \La'_{i_2} + (2-i_2) \\
\vdots & \vdots & \vdots \\
\tau_{n-k} & \equiv & \La'_{i_{n-k}} + (n-k-i_{n-k}). \\
\end{eqnarray*}
Now $\La'_i-i \neq \La'_j -j$ for $i \neq j$.  Furthermore, 
for all $s$ and $t$, 
$| i_s - i_t | \leq n-k-1$, while $|\La'_s - \La'_t | \leq k$ by Remark \ref{rem:ldm}(iv).
Combining these observations, we see that, for $i \neq j$, we have 
$\La'_i-i \not\equiv \La'_j -j \pmod{n}$.  
Therefore, the value of $\tau_j -j$ 
determines $\La'_{i_j} - i_j$, and hence determines $i_j$ for all $j$.  Therefore, $\tau$
determines $r$, implying the result.  

(b) If $\order{\La'+rn}$ is not a decreasing sequence, then $\tau_j +1 = \tau_{j+1}$ for some
$j$.  The congruences above therefore imply that
\[
\La'_{i_j} + (j-i_j) +1 \equiv \La'_{i_{j+1}} + (j +1 - i_{j+1}) \pmod{n}.
\]
which we saw was impossible for $i_j \neq i_{j+1}$.  

(c) Suppose $r_1 \geq 0$.  Then $\order{\La'+rn}_1 \geq \La'_1 + r_1 n \geq \La'_1 > k$, since
$d \geq 1$.  Now suppose $r_1 < 0$.  Therefore, $r_i > 0$ for some $2 \leq i \leq n-k$.  
We have 
\[
\order{\La'+rn}_1 \geq \La'_i + r_i n - (i-1) \geq 0 + n - (n-k-1) = k+1.
\]

(d) Observe that there exists an integer sequence $\rb = (\rb_1, \ldots, \rb_{n-k})$ such that $\order{\increm{\order{\La'+r n}}{1}{-n}} = \order{\La'+\rb n}$.  Thus, (b) implies that
$\order{\increm{\order{\La'+r n}}{1}{-n}}$ is a weakly
decreasing sequence.  It remains to show
that $\order{\increm{\order{\La'+r n}}{1}{-n}}_{n-k} \geq 0$.  Now 
$\order{\increm{\order{\La'+r n}}{1}{-n}}$ is obtained from $\increm{\order{\La'+r n}}{1}{-n}$ 
by applying adjacent transpositions to move the first entry $\order{\La'+rn}_1 -n$ to the right 
until it has no more larger entries to its right.  There are two possibilities.  Either it gets moved
all the way to the $(n-k)$th position, in which case
\[
\order{\increm{\order{\La'+r n}}{1}{-n}}_{n-k} = \order{\La'+rn}_1 -n +(n-k-1) \geq 0
\]
by (c).  Alternatively, it has no larger entries to its right before it 
reaches the $(n-k)$th position, in which case
\[
\order{\increm{\order{\La'+r n}}{1}{-n}}_{n-k} = (\increm{\order{\La'+r n}}{1}{-n})_{n-k}
= \order{\La'+r n}_{n-k} \geq 0
\]
since $\order{\La'+r n}$ is a partition.

\end{proof}

\begin{proof}[Proof of Theorem \ref{thm:gkribbons}]
We use Theorem \ref{thm:gk} as our starting point.
We must show that for any $\ldm$ that is a subposet of $\cnk$, 
\begin{equation}\label{equ:twoformulations}
\sum_{\genfrac{}{}{0pt}{}{r_1 + \cdots +r_{n-k}=0}{r_i \in \mbz}}
s_{(\modif{\La'}{r})'/\mu}(x) = \sum_\tau \vare(\tau/\la) s_{\tau/\mu} (x),
\end{equation}
where $\La = \addribbon{\la}{d}$ and where the sum on the right-hand side  
is over all $\tau$ with $\tau_1 \leq n-k$ that can be obtained from $\la$ by
adding $d$ $n$-ribbons.  First, notice that $\mu$ plays a very
straightforward role.  In particular, if \eqref{equ:twoformulations} holds
for $\dempty{\la}$, then it holds for $\ldm$.  Therefore, we will assume
that $\mu=\emptyset$.  Define two multisets $\Le_d$ and $\R_d$ 
of signed integer sequences as follows:
\begin{eqnarray*}
\Le_d  & = &  \left\{\order{\modif{\addribbon{\la}{d}'}{r}} \in \SPar\ |\ r_i \in \mbz, r_1 + \cdots r_{n-k} = 0 \right\},\\
\R_d & = & \left\{ \vare(\tau/\la) \tau' \in \SPar \ |\ \tau_1 \leq n-k, \mbox{$\tau$ can be obtained}  \right. \\
& & \left. \ \mbox{from $\la$ by adding $d$ $n$-ribbons}\right\}.
\end{eqnarray*}
We see that showing \eqref{equ:twoformulations} amounts to showing that $\Le_d=\R_d$.
Every element of $\R_d$ occurs with multiplicity 1 by definition, and every element of $\Le_d$
occurs with multiplicity 1 by Lemma \ref{lem:mult1}(a).  
So $\Le_d$ and $\R_d$ are, in fact, just sets.

Suppose first that $d=0$, in which case $\R_d = \{\la'\}$.  
Because $\La'=\la'$ and $\la'_i \leq k$ for all $i$, 
we see that the signed sequence $\order{\modif{\La'}{r}}$ will have a negative entry
unless $r_i=0$ for all $i$.  Therefore, $\Le_0=\{\la'\}$ also.

Now suppose that $d > 0$, and assume by induction that $\Le_{d-1}=\R_{d-1}$.  
Define $\IL_d$ to be the set of signed partitions given by
\[
\IL_d = \left\{ (-1)^{n-k-1}\order{\increm{\al}{i}{n}} \in \SPar\ |\ \al \in \Le_{d-1}, 1 \leq i \leq n-k \right\}.
\]
$\IL_d$ can be thought of as an inductive version of $\Le_d$.  
That $\IL_d = \R_d$ is exactly the content of Lemma \ref{lem:simtoribbons}(a), combined
with the induction hypothesis.
It remains to show that $\IL_d = \Le_d$.  

We know that $\Lam = \addribbon{\la}{d-1}$ is obtained from 
$\La = \addribbon{\la}{d}$ by removing an 
\linebreak
$n$-ribbon whose leftmost column is column 
1 and whose rightmost column is column $n-k$.  Lemma \ref{lem:simtoribbons}(b) then
implies that $\Lam' = (-1)^{n-k-1} \order{\increm{\La'}{1}{-n}}$
as signed partitions.  
Also,
Lemma \ref{lem:simtoribbons}(a) implies that
$\La' = (-1)^{n-k-1} \order{\increm{\Lam'}{n-k}{n}}$.

Towards showing that $\IL_d \subseteq \Le_d$, we next consider
$\al = \order{\Lam' + \rb n} \in \Le_{d-1}$ with $\rb_1 + \cdots + \rb_{n-k} = 0$.  
Since 
$\La'_1 \geq \La'_2 \geq \cdots \geq \La'_{n-k} \geq \La'_1 -k$, we get
\begin{eqnarray*}
& & \order{\Lam' + \rb n} \\
& = & (-1)^{n-k-1} \order{\order{\increm{\La'}{1}{-n}} + \rb n} \\
& = & (-1)^{n-k-1} \order{(\La'_2-1, \La'_3-1, \ldots, \La'_{n-k}-1, \La'_1-n+(n-k-1)) + \rb n} \\
& = & (-1)^{n-k-1}
\order{(\La'_2 -1+\rb_1 n, \ldots, \La'_{n-k}-1+\rb_{n-k-1}n,  \La'_1-k-1+\rb_{n-k}n)}  \\
& = & (-1)^{n-k-1} \order{\La' + (\rb_{n-k}-1, \rb_1, \rb_2, \ldots, \rb_{n-k-1})n}.
\end{eqnarray*}
Now suppose we take $\be = (-1)^{n-k-1}\order{\increm{\al}{i}{n}} \in \IL_d$. 
We have
\begin{eqnarray*}
\be & = &  (-1)^{n-k-1}\order{
\increm{(-1)^{n-k-1} \order{\La' + (\rb_{n-k}-1, \rb_1, \rb_2, \ldots, \rb_{n-k-1})n}}{i}{n}} \\
& = & (-1)^{n-k-1}(-1)^{n-k-1} \order{\La' + rn}.
\end{eqnarray*}
for suitable choice of $r=(r_1, \ldots, r_{n-k})$ with $r_1 + \cdots r_{n-k}=0$. 
Therefore, $\be \in \Le_d$ and so
$\IL_d \subseteq \Le_d$.  

Now suppose we take any $\be = \order{\La' +rn} \in \Le_d$.  By Lemma \ref{lem:mult1}(d), 
\linebreak
$\order{\increm{\be}{1}{-n}}=\order{\increm{\order{\La' +rn}}{1}{-n}}$ is a signed partition,
which we choose to denote by $(-1)^{n-k-1}\al$.  It follows that 
$\be = \order{(-1)^{n-k-1}\increm{\al}{i}{n}} = (-1)^{n-k-1}\order{\increm{\al}{i}{n}}$ for some $i$.
Now 
\begin{eqnarray*}
\al & = & (-1)^{n-k-1} \order{\increm{\order{\La' +rn}}{1}{-n}} \\
& = & (-1)^{n-k-1} 
\order{\increm{\order{(-1)^{n-k-1} \order{\increm{\Lam'}{n-k}{n}} +rn}}{1}{-n}} \\
& = & (-1)^{n-k-1}(-1)^{n-k-1}\order{\Lam'+\rb n}  
\end{eqnarray*}
for suitable choice of $\rb = (\rb_1, \ldots, \rb_{n-k})$ with $\rb_1 + \cdots + \rb_{n-k}=0$.  
We conclude that $\al \in \Le_{d-1}$ and hence $\be \in \IL_{d}$.  Therefore
$\Le_{d} \subseteq \IL_{d}$ and so $\Le_d = \IL_d$.  
\end{proof}


\section{Cylindric Schur-positivity}\label{sec:cylschurpos}

Before presenting the conjecture which is the main subject of this
section, we begin with a relevant application of Theorem \ref{thm:gkribbons}.

In the same way that Schur functions are those skew Schur functions $s_{\la/\mu}(x)$
with $\mu=\emptyset$, we will say that \emph{cylindric Schur functions}
are those cylindric skew Schur functions $\cssf{\ldm}(x)$ with $\mu=\emptyset$.
While the Schur functions are known to be a basis for the symmetric functions, we have
the following result for the cylindric Schur functions.

\begin{proposition}
For a given $k, n-k$, the cylindric Schur functions of the form $\cssf{\dempty{\la}}(x)$,
with $\dempty{\la}$ a subposet of $\cnk$, are linearly independent.
\end{proposition}

\begin{proof}
Consider the expansion \eqref{equ:gkribbons} of a cylindric skew Schur function 
$\cssf{\ldm}(x)$ in terms of skew Schur functions.  When $\mu=\emptyset$, 
this expansion is in terms of Schur functions.  Furthermore, a Schur function 
$s_\nu$ can only appear in the Schur expansion of $\cssf{\dempty{\la}}(x)$ if 
$\la$ is the $n$-core of $\nu$.  It follows that when we take a linear
combination of cylindric Schur functions of the form $\cssf{\dempty{\la}}(x)$
having $\dempty{\la}$ a subposet of $\cnk$, we don't get any cancellation 
among the Schur expansions of the cylindric Schur functions.  In particular,
the cylindric Schur functions are linearly independent.  
\end{proof}

We might next ask if every 
cylindric skew Schur function $\cssf{\ldm}(x)$ with $\ldm$ a subposet
of $\cnk$ can be expressed 
as a linear combination of cylindric Schur functions of the form
$\cssf{\nu/e/\emptyset}(x)$, where each 
$\nu/e/\emptyset$ 
is also a subposet of $\cnk$.
As we shall see, an
affirmative answer to this question would also imply Conjecture
\ref{con:cylschurpos} below.  

\begin{definition}
Suppose $\ldm$ is a cylindric skew shape that is a subposet of $\cnk$.
We say that $\cssf{\ldm}(x)$ in the variables $x=(x_1, x_2, \ldots)$ is
\emph{cylindric Schur-positive} if it can be expressed as a linear 
combination of cylindric Schur functions
$\cssf{\nu/e/\emptyset}(x)$ with positive coefficients, 
where each such $\nu/e/\emptyset$ is also a subposet of $\cnk$.  
\end{definition}

As an analogue of the fact that every skew Schur function is Schur-positive, 
we propose the following conjecture.

\begin{conjecture}\label{con:cylschurpos}
Every cylindric skew Schur function is cylindric 
Schur-positive. 
\end{conjecture}

As we noted in Remark \ref{rem:ribbonspositive}, this conjecture is true for
cylindric ribbons.  The rest of this section will be devoted to other evidence in favour of
the conjecture.  

It follows from \eqref{equ:partsizes} that we can split $\cssf{\ldm}(x)$ into two 
sums as follows:
\begin{equation}\label{equ:twosums}
\cssf{\ldm}(x) = \sum_{\nu \inpkn} a_\nu s_\nu(x) + \sum_{\genfrac{}{}{0pt}{}{\nu:\nu_1\leq n-k}{l(\nu)>k}}
b_\nu s_\nu(x).
\end{equation}
When $\nu \inpkn$, we know that $s_\nu(x)$ is a cylindric Schur function.  Furthermore, we know
from Theorem \ref{thm:postnikov} that $a_\nu \geq 0$ for all $\nu \inpkn$.  Therefore, the
first sum is cylindric Schur-positive.  

Now consider the second sum, which we denote by $\error{\ldm}{x}$.
We know that $\cssf{\ldm}(x)$ is cylindric Schur-positive
when $d=0$.  Therefore, we can assume by induction that $\cssf{\la/(d-1)/\mu}(x)$ is
cylindric Schur-positive:
\begin{equation}\label{equ:induction}
\cssf{\la/(d-1)/\mu}(x) = \sum_{\genfrac{}{}{0pt}{}{\nu,e}{\nu \inpkn}} c_{\nu,e} \cssf{\nu/e/\emptyset}(x),
\end{equation}
where $c_{\nu,e} \geq 0$ for all $\nu, e$, and $e$ is a always non-negative integer.
(For $\cssf{\nu/e/\emptyset}(x) \neq 0$, we require that
$ne=|\la|-|\mu|+n(d-1)-|\nu|$.)
We conjecture, in fact, that $\error{\ldm}{x}$ can be expressed exactly in terms of 
$\cssf{\la/(d-1)/\mu}(x)$ as:
\[
\error{\ldm}{x} = \sum_{\genfrac{}{}{0pt}{}{\nu,e}{\nu \inpkn}} c_{\nu,e} \cssf{\nu/e+1/\emptyset}(x).
\]
Plugging this into \eqref{equ:twosums}, we get
\begin{equation}\label{equ:errorterm}
\cssf{\ldm}(x) = \sum_{\nu \inpkn} a_\nu s_\nu(x) +
\sum_{\genfrac{}{}{0pt}{}{\nu,e}{\nu \inpkn}} c_{\nu,e} \cssf{\nu/e+1/\emptyset}(x),
\end{equation}
where $a_\nu, c_{\nu, e} \geq 0$ for all $\nu, e$.  
This expression is a strong refinement of Conjecture \ref{con:cylschurpos} as it 
gives much information about the form  
of the cylindric 
Schur-positive expansion of $\cssf{\ldm}(x)$.  
Using \cite{BucSoftware,SteSoftware} we have verified \eqref{equ:errorterm} for all $\ldm$ with 
$k, n-k, d \leq 5$.  

One way to show \eqref{equ:errorterm} would be to show that the coefficient of $s_\si(x)$ 
is the same on both sides for all partitions $\si$ with $|\si| = |\la| + nd -|\mu|$. 
Since we are only worried about the cylindric Schur-positivity of $\error{\ldm}{x}$, assume
that $\si_1 \leq n-k$ but $l(\si) > k$.  There is a certain important class of such 
partitions $\si$
for which we can show $s_\si(x)$ has the same coefficient on both sides of 
\eqref{equ:errorterm}:

\begin{proposition}\label{pro:somecoeffs}
Suppose we are given a cylindric shape $\ldm$ which is a subposet of $\cnk$ and,
to avoid trivialities, we take $d \geq 1$.   Consider a partition $\si$ with
$|\si| = |\la| + nd -|\mu|$, $\si_1 \leq n-k$, $l(\si) > k$ and the additional condition
that 
\begin{equation}\label{equ:cylindricrim}
\si'_1 \geq \si'_2 \geq \cdots \geq \si'_{n-k} \geq \si'_1 -k.
\end{equation}
Then 
\[
[s_\si(x)]\cssf{\ldm}(x) = [s_\si(x)]
\sum_{\genfrac{}{}{0pt}{}{\nu,e}{\nu \inpkn}} c_{\nu,e} \cssf{\nu/e+1/\emptyset}(x),
\]
where 
\[
\cssf{\la/(d-1)/\mu}(x) = \sum_{\genfrac{}{}{0pt}{}{\nu,e}{\nu \inpkn}} c_{\nu,e} \cssf{\nu/e/\emptyset}(x).
\]
\end{proposition}

\begin{proof}
The key idea is that since $\si$ satisfies \eqref{equ:cylindricrim}, $\la/d/\si$ is a valid 
cylindric shape.  Because of the conditions on $\si$, we also know that 
$\addribbon{\si}{-1}$ is a well-defined partition.  Indeed, by \eqref{equ:cylindricrim} we
know that $\si'_{n-k} > 0$ and
\[
\addribbon{\si}{-1}' = (\si'_2-1, \ldots, \si'_{n-k}-1, \si'_1-k-1).  
\]
By Theorem \ref{thm:gkribbons}, we have
\begin{eqnarray*}
[s_\si(x)]\cssf{\ldm}(x) & =  & [s_\si(x)]\sum_\tau \vare(\tau/\la) s_{\tau/\mu} (x) \\
& = & \sum_\tau \vare(\tau/\la) c^\tau_{\mu \si} \\
& = & [s_\mu(x)]\sum_\tau \vare(\tau/\la) s_{\tau/\si} (x) \\
& = & [s_\mu(x)]\cssf{\la/d/\si}(x),
\end{eqnarray*}
where the sums are over all $\tau$ with $\tau_1 \leq n-k$ that can be obtained from $\la$ by
adding $d$ $n$-ribbons, and where $c^\tau_{\mu \si}$ denotes the Littlewood-Richardson
coefficient.  By Remark \ref{rem:ldm}(ii), $\la/d/\si$ and $\la/(d-1)/\addribbon{\si}{-1}$ are
the same cylindric skew shape.  Therefore, now with the sums over 
all $\tau$ with $\tau_1 \leq n-k$ that can be obtained from $\la$ by
adding $d-1$ $n$-ribbons, we have
\begin{eqnarray*}
[s_\si(x)]\cssf{\ldm}(x) & =  & [s_\mu(x)]\cssf{\la/(d-1)/\addribbon{\si}{-1}}(x) \\
& = & [s_\mu(x)]\sum_\tau \vare(\tau/\la) s_{\tau/\addribbon{\si}{-1}} (x) \\
& = & \sum_\tau \vare(\tau/\la) c^\tau_{\mu \addribbon{\si}{-1}} \\
& = &  [s_{\addribbon{\si}{-1}}(x)] \cssf{\la/(d-1)/\mu}(x) \\
& = & [s_{\addribbon{\si}{-1}}(x)]
\sum_{\genfrac{}{}{0pt}{}{\nu,e}{\nu \inpkn}} c_{\nu,e} \cssf{\nu/e/\emptyset}(x).
\end{eqnarray*}
However, since $\si$ and $\addribbon{\si}{-1}$ have the same $n$-core, 
\[
[s_{\addribbon{\si}{-1}}(x)]
\sum_{\genfrac{}{}{0pt}{}{\nu,e}{\nu \inpkn}} c_{\nu,e} \cssf{\nu/e/\emptyset}(x) = 
[s_\si(x)] \sum_{\genfrac{}{}{0pt}{}{\nu,e}{\nu \inpkn}} c_{\nu,e} \cssf{\nu/e+1/\emptyset}(x),
\]
as required.
\end{proof}

As promised, we can now reformulate Conjecture \ref{con:cylschurpos} into 
a seemingly easier statement.

\begin{corollary}
Conjecture \ref{con:cylschurpos} holds if and only if every 
cylindric skew Schur function $\cssf{\ldm}(x)$ with $\ldm$ a subposet
of $\cnk$ can be expressed 
as a linear combination of cylindric Schur functions
$\cssf{\nu/e/\emptyset}(x)$, where each 
$\nu/e/\emptyset$ 
is also a subposet of $\cnk$.
\end{corollary}

In other words, to prove Conjecture \ref{con:cylschurpos}, 
we don't have to show that the coefficients are positive.

\begin{proof}
The ``only if'' direction is trivial.  So suppose $\cssf{\ldm}(x)$ can be expressed as a linear
combination of cylindric Schur functions.    Let $\cssf{\nu/m/\emptyset}(x)$ be a cylindric
Schur function that appears with coefficient $a_{\nu, m}$ in this linear combination.  
We need to show that $a_{\nu, m} \geq 0$.  Assume that 
$|\la| + dn -|\mu| = |\nu|+mn$, since otherwise $a_{\nu, m} = 0$.  We proceed by induction
on $d$, with the case $d=0$ being trivial.  

If $m=0$, we know by \eqref{equ:twosums} that $a_{\nu, m} \geq 0$.  Therefore, assume that
$m \geq 1$.
Consider $\si = \addribbon{\nu}{m}$.  
Using the fact that $\nu \inpkn$, we can check that $\si$ satisfies
\[
\si'_1 \geq \si'_2 \geq \cdots \geq \si'_{n-k} \geq \si'_1 -k.
\]
Therefore, we can apply Proposition \ref{pro:somecoeffs}.
We get that 
\begin{equation}\label{equ:somecoeffs}
[s_\si(x)]\cssf{\ldm}(x) =
[s_\si(x)]
\sum_{\genfrac{}{}{0pt}{}{\tau,e}{\tau \inpkn}} c_{\tau,e} \cssf{\tau/e+1/\emptyset}(x),
\end{equation}
where 
\[
\cssf{\la/(d-1)/\mu}(x) = \sum_{\genfrac{}{}{0pt}{}{\tau,e}{\tau \inpkn}} c_{\tau,e} \cssf{\tau/e/\emptyset}(x).
\]
By the induction hypothesis, $c_{\tau, e} \geq 0$ for all $\tau, e$.  
Since $\si$ has $n$-core $\nu$, we know from Theorem \ref{thm:gkribbons}
that $s_\si(x)$ appears with 
coefficient $\vare(\si/\nu)=1$ in $\cssf{\nu/m/\emptyset}(x)$ and appears
with coefficient 0 in the Schur expansion of any other cylindric Schur function. 
Therefore, \eqref{equ:somecoeffs} tells us that
$a_{\nu, m} = c_{\nu, m-1} \geq 0$, as required.

\end{proof}
 

\bibliography{../../master}

\newcommand{\noopsort}[1]{} \newcommand{\printfirst}[2]{#1}
  \newcommand{\singleletter}[1]{#1} \newcommand{\switchargs}[2]{#2#1}
\begin{thebibliography}{10}

\bibitem{BeKn72}
Edward~A. Bender and Donald~E. Knuth.
\newblock Enumeration of plane partitions.
\newblock {\em J. Combinatorial Theory Ser. A}, 13:40--54, 1972.

\bibitem{BCF99}
Aaron Bertram, Ionu{\c{t}} Ciocan-Fontanine, and William Fulton.
\newblock Quantum multiplication of {S}chur polynomials.
\newblock {\em J. Algebra}, 219(2):728--746, 1999.

\bibitem{BucSoftware}
Anders~S. Buch.
\newblock {L}ittlewood-{R}ichardson calculator, 1999.
\newblock \newline Available from \texttt{http://home.imf.au.dk/abuch/lrcalc/}.

\bibitem{Cum91}
C.~J. Cummins.
\newblock {${\rm su}(n)$} and {${\rm sp}(2n)$} {WZW} fusion rules.
\newblock {\em J. Phys. A}, 24(2):391--400, 1991.

\bibitem{Gei77}
Ladnor Geissinger.
\newblock Hopf algebras of symmetric functions and class functions.
\newblock In {\em Combinatoire et repr\'esentation du groupe sym\'etrique
  (Actes Table Ronde C.N.R.S., Univ. Louis-Pasteur Strasbourg, Strasbourg,
  1976)}, pages 168--181. Lecture Notes in Math., Vol. 579. Springer, Berlin,
  1977.

\bibitem{Ges84}
Ira~M. Gessel.
\newblock Multipartite {$P$}-partitions and inner products of skew {S}chur
  functions.
\newblock In {\em Combinatorics and algebra (Boulder, Colo., 1983)}, volume~34
  of {\em Contemp. Math.}, pages 289--317. Amer. Math. Soc., Providence, RI,
  1984.

\bibitem{GeKr97}
Ira~M. Gessel and Christian Krattenthaler.
\newblock Cylindric partitions.
\newblock {\em Trans. Amer. Math. Soc.}, 349(2):429--479, 1997.

\bibitem{GoWe90}
Frederick~M. Goodman and Hans Wenzl.
\newblock Littlewood-{R}ichardson coefficients for {H}ecke algebras at roots of
  unity.
\newblock {\em Adv. Math.}, 82(2):244--265, 1990.

\bibitem{Kac90}
Victor~G. Kac.
\newblock {\em Infinite-dimensional {L}ie algebras}.
\newblock Cambridge University Press, Cambridge, third edition, 1990.

\bibitem{Mac95}
I.~G. Macdonald.
\newblock {\em Symmetric functions and {H}all polynomials}.
\newblock Oxford Mathematical Monographs. The Clarendon Press Oxford University
  Press, New York, second edition, 1995.

\bibitem{Mal93}
Claudia Malvenuto.
\newblock {$P$}-partitions and the plactic congruence.
\newblock {\em Graphs Combin.}, 9(1):63--73, 1993.

\bibitem{MalThesis}
Claudia Malvenuto.
\newblock {\em Produits et coproduits des fonctions quasi-sym\'etriques et de
  l'alg\`ebre des descentes}, volume~16 of {\em Publications du Laboratoire de
  Combinatoire et d'Informatique Math\'ematique}.
\newblock Laboratoire de Combinatoire et d'Informatique Math\'ematique,
  Universit\'e du Qu\'ebec \`a Montr\'eal, 1993.
\newblock Ph.D. thesis.

\bibitem{MaRe95}
Claudia Malvenuto and Christophe Reutenauer.
\newblock Duality between quasi-symmetric functions and the {S}olomon descent
  algebra.
\newblock {\em J. Algebra}, 177(3):967--982, 1995.

\bibitem{McNThesis}
Peter McNamara.
\newblock {\em Edge labellings of partially ordered sets}.
\newblock PhD thesis, Massachusetts Institute of Technology, 2003.
\newblock \texttt{http://www.lacim.uqam.ca/\symbol{126}mcnamara/research.html}.

\bibitem{Pos04pr}
Alexander Postnikov.
\newblock Affine approach to quantum {S}chubert calculus.
\newblock {\em Duke Math. J.}, to appear.
\newblock \texttt{http://www.arxiv.org/abs/math.CO/0205165}.

\bibitem{StaThesis}
Richard~P. Stanley.
\newblock {\em Ordered structures and partitions}.
\newblock American Mathematical Society, Providence, R.I., 1972.
\newblock Memoirs of the American Mathematical Society, No. 119.

\bibitem{ec1}
Richard~P. Stanley.
\newblock {\em Enumerative combinatorics. {V}ol. {I}}.
\newblock Wadsworth \& Brooks/Cole Advanced Books \& Software, Monterey, CA,
  {\noopsort{1997a}}1986.
\newblock Second printing, Cambridge University Press, Cambridge/New York,
  1997.

\bibitem{ec1errata}
Richard~P. Stanley.
\newblock Errata and addenda to \textit{{E}numerative combinatorics. {V}ol. 1,
  second printing}.
\newblock \texttt{http://www-math.mit.edu/\symbol{126}rstan/ec/newerr.ps},
  {\noopsort{1997b}}version of 9th April 2003 or later.

\bibitem{ec2}
Richard~P. Stanley.
\newblock {\em Enumerative combinatorics. {V}ol. 2}, volume~62 of {\em
  Cambridge Studies in Advanced Mathematics}.
\newblock Cambridge University Press, Cambridge, 1999.

\bibitem{Sta02}
Richard~P. Stanley.
\newblock Recent developments in algebraic combinatorics.
\newblock {\em Israel J. Math.}, 143:317--339, 2004.

\bibitem{SteSoftware}
John~R. Stembridge.
\newblock {SF}, posets and coxeter/weyl.
\newblock \newline Available from
  \texttt{http://www.math.lsa.umich.edu/\symbol{126}jrs/maple.html}.

\bibitem{Thi91}
Jean-Yves Thibon.
\newblock Coproduits de fonctions sym\'etriques.
\newblock {\em C. R. Acad. Sci. Paris S\'er. I Math.}, 312(8):553--556, 1991.

\bibitem{Wal90}
Mark~A. Walton.
\newblock Fusion rules in {W}ess-{Z}umino-{W}itten models.
\newblock {\em Nuclear Phys. B}, 340(2-3):777--790, 1990.

\bibitem{Zel81}
Andrey~V. Zelevinsky.
\newblock {\em Representations of finite classical groups}, volume 869 of {\em
  Lecture Notes in Mathematics}.
\newblock Springer-Verlag, Berlin, 1981.
\newblock A Hopf algebra approach.

\end{thebibliography}
\bibliographystyle{plain}

\end{document}